\documentclass[12pt,a4paper]{article}
\usepackage[left=1.2in,right=1.2in,top=.8in,bottom=.8in]{geometry}
\usepackage{amsmath,amsfonts,amsthm,amssymb,graphicx,xcolor}  

\newtheorem*{theorem*}{Theorem}
\newtheorem*{proposition*}{Proposition}
\newtheorem{lemma}{Lemma}

\theoremstyle{remark}
\newtheorem{remark}{Remark}
\usepackage[utf8]{inputenc}
\usepackage{color}
\usepackage{hyperref}
\usepackage[mathscr]{euscript}
\usepackage{enumerate}
\usepackage{multirow}

\usepackage{mathtools}
\mathtoolsset{showonlyrefs=true} % show or hide unreferenced

\newcommand{\e}{\mathrm{e}}

\renewcommand{\d}{\mathrm{d}}

\newcommand{\R}{\mathbb{R}}

\newcommand{\E}{\mathcal{E}}

\newcommand{\eps}{\varepsilon}

\newcommand{\N}{\mathbb{N}}

\newcommand{\s}{\hspace{0.5pt}}

\newcommand{\tildeR}{\widetilde{\mathcal{R}}}

\newcommand{\y}{\mathbf{y}}
\newcommand{\g}{\mathbf{g}}

\newcommand{\p}{\partial}

\newcommand{\norm}[1]{\lVert #1 \rVert}

\usepackage{authblk}
\begin{document}
\title{An inverse problem for a semi-linear wave equation: a numerical study}
\author[1]{Matti Lassas}
%%\address{Department of Mathematics and Statistics, University of Jyv\"askyl\"a, Jyv\"askyl\"a, Finland}
%%%\curraddr{}
%%\email{matti.lassas@helsinki.fi}
%
%
\author[2]{Tony Liimatainen}
%\address{Department of Mathematics and Statistics, University of Jyv\"askyl\"a, Jyv\"askyl\"a, Finland}
%%%\curraddr{}
%%\email{tony.liimatainen@helsinki.fi}
%
\author[3]{Leyter Potenciano-Machado}
%%\address{Department of Mathematics and Statistics, University of Jyv\"askyl\"a, Jyv\"askyl\"a, Finland}
%%\email{leyter.m.potenciano@jyu.fi}

\author[4]{Teemu Tyni}
%\address{Department of Mathematics and Statistics, University of Jyv\"askyl\"a, Jyv\"askyl\"a, Finland}
%\email{teemu.tyni@helsinki.fi}

\affil[1]{Department of Mathematics and Statistics, University of Helsinki, Helsinki, Finland}
\affil[2,3]{Department of Mathematics and Statistics, University of Jyv\"askyl\"a, Jyv\"askyl\"a, Finland}
\affil[4]{Department of Mathematics, University of Toronto, Toronto, Canada}
%%\email{}
%\affil[1]{matti.lassas@helsinki.fi}
%\affil[2]{tony.liimatainen@helsinki.fi}
%\affil[3]{leyter.m.potenciano@jyu.fi}
%\affil[4]{teemu.tyni@helsinki.fi}

%\curraddr{}

%\date{} %add this to remove the date from submission version

\maketitle

\begin{abstract}
We consider an inverse problem of recovering a potential associated to a semi-linear wave equation with a quadratic nonlinearity in $1+1$ dimensions.
% 
% % on Rn+1, n≥1. We show a Hölder stability estimate for the recovery of an unknown potential a of the wave equation □u+aum=0 from its Dirichlet-to-Neumann map. We show that an unknown potential a(x,t), supported in Ω×[t1,t2], of the wave equation □u+aum=0 can be recovered in a Hölder stable way from the map u|∂Ω×[0,T]↦⟨ψ,∂νu|∂Ω×[0,T]⟩L2(∂Ω×[0,T]). This data is equivalent to the inner product of the Dirichlet-to-Neumann map with a measurement function ψ. We also prove similar stability result for the recovery of a when there is noise added to the boundary data. The method we use is constructive and it is based on the higher order linearization. As a consequence, we also get a uniqueness result. We also give a detailed presentation of the forward problem for the equation □u+aum=0. 
% 
% 
% We consider an inverse problem for one-dimensional semilinear wave equation with a quadratic nonlinearity.
%
We develop a numerical scheme to determine the potential from a noisy Dirichlet-to-Neumann map on the lateral boundary.  The scheme is based on the recent higher order linearization method \cite{LLPT20}. 
%incorporated with noise.
%
We also present an approach to numerically estimating two-dimensional derivatives of noisy data via Tikhonov regularization.
The methods are tested using synthetic noisy measurements of the Dirichlet-to-Neumann map. Various examples of reconstructions of the potential functions are given.
\end{abstract}
This paper is a continuation of \cite{LLPT20} into a computational direction. We study an inverse boundary value problem for the one-dimensional non-linear wave equation
\begin{equation}\label{eq:intro_wave-eq}
\begin{cases}
\frac{\partial^2}{\partial t^2}u(x,t) = \frac{\partial}{\partial x^2}u(x,t) - q(x,t)\s u(x,t)^{2}  &(x,t)\in [a,b]\times[0,T],\\
u(a,t)=f^L(t), &t\in [0,T], \\
u(b,t) = f^R(t), &t\in [0,T], \\
u(x,0) = \partial_t u(x,0) = 0 &x\in [a,b] .
\end{cases}
\end{equation}
It was shown in \cite{LLPT20} that there exists a unique small solution to \eqref{eq:intro_wave-eq} if $f^L,f^R\in H^s([0,T])$ have small enough norms in $H^s([0,T])$ for sufficiently large $s>0$.
%
%The existence of a unique solution, for small enough data $f^L,f^R\in H^s([0,T])$, can be found from \cite{LLPT20}.
%
Let us denote the lateral boundary of $[a,b]\times [0,T]$ by 
\[
\Sigma := \{a,b\}\times [0,T].
\]
The Sobolev space $H^{s+1}(\Sigma)$ is naturally isomorphic to $H^{s+1}([0,T])\times H^{s+1}([0,T])$.
The Dirichlet-to-Neumann map (DN map) $\Lambda_q$ is then defined as
\begin{align*}
\Lambda_q &: H^{s+1}([0,T])\times H^{s+1}([0,T]) \to H^s(\Sigma),\\
 &\quad\{f^L,f^R\} \mapsto \p_\nu u\big|_\Sigma,
\end{align*}
where $u$ solves \eqref{eq:intro_wave-eq} and $\nu$ is the outward pointing normal of $\Sigma$. The notations $f^L$ and $f^R$ refer to the boundary value on the left and right side of $\Sigma$ as in \eqref{eq:intro_wave-eq}. % and $N$ is the outward pointing normal of $\Sigma$.
The inverse problem we consider in this paper is the recovery of the unknown potential function $q\in C^\infty([a,b]\times [0,T])$ from the DN map $\Lambda_q$.
It was shown in \cite{LLPT20} that the DN map uniquely determines the unknown potential $q(x,t)$. The work also provided a stable reconstruction algorithm. In this work we implement the reconstruction algorithm numerically. To the best of our knowledge, this is the first work that provides numerical results based on using non-linearity as a beneficial tool for inverse problems of nonlinear wave equations. % equation \eqref{eq:intro_wave-eq}.  

Our numerical reconstruction algorithm follows the recent work \cite{LLPT20}. The work is based on the recent \emph{higher order linearization method}, which uses boundary values with several parameters (two in this paper) and obtains new linearized equations after differentiating with respect to these parameters. The reconstruction algorithm based on the method is fast and straightforward to implement. Given a suitable measurement data set, the inversion is done practically in real time. We explain the higher order linearization method and the the numerical reconstruction in Section \ref{sec:prelim}. 
%To the best of our knowledge this work is the first one that uses the higher order linearization method for numerics of inverse problems of a boundary value problem for nonlinear wave equations. 
Our  reconstruction works equally well for both time independent and time dependent potentials. 
%~\f{Let's discuss the benefits in short here. Let's give couple of references in the linear case where the time dependent potential is difficult. Matti, do you know?}

%By using the nonlinearity as a tool, some still unsolved inverse problems for hyperbolic linear equations
%have been solved for their nonlinear counterparts. 

 %The idea of a first or a second order linearization was initiated in \cite{isakov1993uniqueness,isakov1994global} to study inverse problems for nonliner equations. 
 The idea of higher order linearizations was developed by Kurylev, Lassas and Uhlmann \cite{KLU18}.
 %in the context of nonlinear hyperbolic equations over Lorentzian geometries.
%The higher order linearization method for nonlinear wave equations equations originates from the work \cite{KLU18} by Kurylev, Lassas and Uhlmann. 
They observed that non-linearity can be used as a beneficial tool in inverse problems and recovered a Lorentzian manifold up to a conformal change from local measurements for the scalar wave equation with quadratic nonlinearity. The measurements in the work were modeled by a source-to-solution mapping, which assigns to a source the corresponding solution of the nonlinear wave equation restricted to the measurement set.

Instead of measurements modeled by source-to-solution mapping, the higher order linearization method evolved to a tool for inverse problems for boundary value problems in the context of elliptic equations in \cite{FO19, KrUh19, KrUh20, LLLS19a, LLLS19b}. Especially the works \cite{KrUh19,  KrUh20, LLLS19b} introduced the concept of \emph{measurement function} (the term measurement function was coined in \cite{LLPT20}) that allowed to apply the higher order linearization method also for inverse problems of boundary value problems for nonlinear wave equations.  
%The research of inverse problems for non-linear equations using the higher order linearization method is expanding fast. 
The works \cite{HUZ20, HUZ21, LLPT21} study inverse problems for boundary value problems for nonlinear wave equations by using the aforementioned method. We refer the reader to the works \cite{BKLT20, dH2019, dH2020, FLO21, FO20, LaUhYa20, LUW17,LUW18, UZ21,UZ21b, WZ2019} for more examples of inverse problems for nonlinear wave type and hyperbolic equations, and to \cite{HUZ21} for additional references. The very recent work \cite{SSX21} does a numerical study of an inverse problem for a nonlinear elliptic equation by using the higher order linearization method.

In the current paper, we use the higher order linearization method to measure the response coming from nonlinear interactions of waves that approximate delta functions. %A possible downside of this approach is that we rely on measurements which are thus point-wise, that is, 
This means that for each $(x_0,t_0)\in [a,b]\times [0,T]$, where we wish to recover the unknown potential $q$, we need to make one measurement. One can think that a possible downside of this approach is that since it relies in this sense point-wise measurements, it does not average out well the noise in the measurements. However, as we will observe, our reconstruction is  good even in the presence of a reasonable noise. We will discuss this matter in Section \ref{sec:IP}.
%
% Heuristically, the reason for this is that we probe the unknown potential by plane waves, which interact in the presence of the nonlinearity at the point $(x_0,t_0)$ and produce what is essentially a virtual delta source at that point.
% %
% Once we measure the signal produced by this interaction at the boundary, we are able to recover the potential function $q$.
%

%
The structure of this paper is as follows. In Section~\ref{sec:prelim} we first discuss the theoretical results of \cite{LLPT20}. %There we also define the boundary values $f^R$ and $f^L$ of the equation \eqref{eq:intro_wave-eq}, which we will use in the numerical reconstruction. 
After that we discuss how the theoretical results are translated into a numerical recovery method. 
Section~\ref{sec:forward} concerns the numerical implementation of solving the forward problem \eqref{eq:intro_wave-eq}. The forward problem is solved by a simple finite difference method and it is needed to build the synthetic DN map.
We present an example and discuss the convergence of the finite difference method.
In Section~\ref{sec:Dregularization} we discuss an approach for approximating two dimensional derivatives by using a regularization method. 
In Section~\ref{sec:IP} we discuss a numerical implementation of the reconstruction algorithm for the recovery of the potential $q$ from the associated DN map.
Finally, in Section~\ref{sec:examples} we present and discuss examples of recoveries of various potentials.
\section{Preliminaries}\label{sec:prelim}

\subsection{The higher order linearization method}
Here we explain how to recover $q(x,t)$ theoretically from the DN map $\Lambda_q$ associated with~\eqref{eq:intro_wave-eq}. As mentioned, a theoretical method for recovering $q$ was established in the authors' earlier work~\cite{LLPT20}, where stability of the recovery was also considered. The recovery in~\cite{LLPT20} was based on \emph{the higher order linearization method}, which we now explain. We follow the notation of the work~\cite{LLPT20} and refer to the work for justifications and details of the following formal discussion. % described in the introduction. ~\f{References somewhere.} 

% The following formal argument is justified theoretically in~\cite{LLPT20}.
%¤already given in The recocery method we use was already given 
%authors earlier~\cite{LLPT20}. 

%We will also use an integration by parts argument introduced in the study of partial data inverse problem for non-linear elliptic equations in \cite{LLLS19b, KrUh20}. Similar argument was also used recently in~\cite{HUZ20}.

%\cite{AsZh17, CaNaVa19, FeOk20, KaNa02, LLLS19a, LLLS19b, KrUh19, KrUh20, SunUh97}.
%

%We first explain how we can recover $q$ from the DN map $\Lambda$ of the equation~\eqref{eq:intro_wave-eq}.
%This method can be used to %by presenting a formal argument how to 
%recover $q$ from the knowledge of the DN map $\Lambda$ of the equation~\eqref{eq:intro_wave-eq} as follows.
 %Let us consider the case $m=2$. 
 
 Let $f_1,\s f_2\in H^{s+1}(\Sigma)$, $s>1$. %, and let us denote by $u_{\eps_1f_1+\eps_2f_2}$ the solution to~\eqref{eq:intro_wave-eq} with boundary data $\eps_1f_1+\eps_2f_2$, where $\eps_1,\eps_2$ are sufficiently small parameters. 
 To avoid confusion with the notation, we recall that $f_1=(f_1^L,f_1^R)$ and  $f_2=(f_2^L,f_2^R)$, where each component function belongs to $H^{s+1}([0,T])$. Let us consider the family $u_{\eps_1f_1+\eps_2f_2}$ of solutions to
 \begin{equation}\label{eq:prelim_wave-eq_recovery}
\begin{cases}
%\frac{\partial^2}{\partial t^2}u_{\eps_1f_1+\eps_2f_2}(x,t) = \frac{\partial}{\partial x^2}u_{\eps_1f_1+\eps_2f_2}(x,t) - q(x,t)\s u_{\eps_1f_1+\eps_2f_2}^{2}(x,t)   = 0 &(x,t)\in [a,b]\times[0,T], \\
\square u_{\eps_1f_1+\eps_2f_2}(x,t) = - q(x,t)\s u_{\eps_1f_1+\eps_2f_2}^{2}(x,t),   &(x,t)\in [a,b]\times[0,T], \\
u_{\eps_1f_1+\eps_2f_2}(x,t)=\eps_1f_1(x,t)+\eps_2f_2(x,t), &(x,t)\in \Sigma, \\
%u(b,t) = f_2^R(t), &t\in [0,T], \\
u_{\eps_1f_1+\eps_2f_2}(x,0) = \partial_t u_{\eps_1f_1+\eps_2f_2}(x,0) = 0, &x\in [a,b],
\end{cases}
\end{equation}
parametrized by small parameters $\eps_1,\eps_2\in \R$. Here we have denoted
\[
 \square := \frac{\partial^2}{\partial t^2}-\frac{\partial}{\partial x^2}.
\]
We differentiate the equation~\eqref{eq:prelim_wave-eq_recovery} with respect to the parameters $\eps_1$ and $\eps_2$. A formal calculation (justified in~\cite{LLPT20}) shows that the mixed derivative
\begin{equation}\label{eq:w}
 w:=\frac{\p}{\p \eps_1}\frac{\p}{\p \eps_2}\Big|_{\eps_1=\eps_2=0} u_{\eps_1 f_1 + \eps_2 f_2}
\end{equation}
is a function whose initial and boundary data vanish and that satisfies
\begin{equation}\label{eq:second_deriv}
 \frac{\partial^2}{\partial t^2}w(x,t) = \frac{\partial}{\partial x^2}w(x,t) -2q(x,t) v_{1}(x,t)v_2(x,t).
\end{equation}
 %By differentiating the data as well, the initial and boundary data of $w$ vanishes. 
 %Here $v_1$ and $v_2$ are solutions to the linearized equation $\square v=0$ with $v_1|_\Sigma=f_1$ and $v_2|_\Sigma=f_2$. 
 The functions $v_j$ in~\eqref{eq:second_deriv}, for $j=1,2$, solve the linear wave equation
 %~\f{Think what's the best notation for the boundary value arguments.}
 \begin{equation}\label{eq:vj_equations}
 \begin{cases}
 \frac{\partial^2}{\partial t^2}v_j(x,t) = \frac{\partial}{\partial x^2}v_j(x,t), &(x,t)\in [a,b]\times [0,T],\\ 
 v_j(x,t) = f_j(x,t), & (x,t)\in \Sigma, \\%\text{on } \p\Omega\times[0,T],\\
 v_j(x,0) = \p v_j(x,0) = 0, & x\in [a,b].
 \end{cases}
 \end{equation}
 %In the literature, taking mixed mixed derivatives of the equationthe solution $u_{\eps_1f_1+\eps_2f_2}$ of the non-linear equation with respect to $\eps_1$ and $\eps_2$ is known as the higher order linearization method.
  This way we have produced new linear equations from the non-linear equation \eqref{eq:intro_wave-eq}. Note that the equation for $v_j$ is independent of the unknown potential $q(x,t)$. The function $w$ defined by~\eqref{eq:w} is called the second linearization of $u_{\eps_1f_1+\eps_2f_2}$.

  If the DN map is known, then the normal derivative $\p_\nu w$ on the lateral boundary $\Sigma$ of the second linearization is also known, because  %Taking the mixed derivative of the DN map yields
  \[
   \p_\nu w = \p_\nu \left(\p^2_{\eps_1\eps_2}|_{\eps_1=\eps_2=0}\s u_{\eps_1f_1+\eps_2f_2}\right)= \p^2_{\eps_1\eps_2}|_{\eps_1=\eps_2=0}\s \Lambda_q(\eps_1f_1+\eps_2f_2).
 \]
 (See \cite{LLPT20} for justification of this calculation.)
  We let $v_0$ be an auxiliary function solving $\square v_0=0$ in $[a,b]\times [0,T]$ with $v_0(x,T)=\p_tv_0(x,T)=0$ for $x\in [a,b]$. The function $v_0$ is called the \emph{measurement function}. The function $v_0$ is used to compensate the fact that $\p_\nu w$ is not known on $\{t=T\}$ from the DN map $H^{s+1}(\Sigma)\to H^s(\Sigma)$. By multiplying~\eqref{eq:second_deriv} by $v_0$ and integrating by parts on $[a,b]\times [0,T]$, we arrive at the integral identity
 \begin{equation}\label{eq:integral_identity_derivitve}
 \begin{split}
  \int_\Sigma v_0\s \p^2_{\eps_1\eps_2}|_{\eps_1=\eps_2=0}\Lambda_q(\eps_1f_1+\eps_2f_2)\s \d \Sigma &=\int_{[a,b]\times [0,T]}\s v_0 \s \square w \s \d x \s \d t \\
  &= - 2\int_{[a,b]\times [0,T]} q\s v_0\s v_1\s v_2\s \d x \s \d t.
  \end{split} 
 \end{equation}
 Thus the quantity 
 \begin{equation}\label{eq:integral_density}
  \int_{[a,b]\times [0,T]} q\s v_0\s v_1\s v_2 \s \d x \s \d t 
 \end{equation}
 is known from the knowledge of the DN map. The integral identity \eqref{eq:integral_identity_derivitve} holds for arbitrary solutions  $v_0, v_1$ and $v_2$ to $\square v=0$  in $[a,b]\times [0,T]$ satisfying the initial conditions $v_0(x,T) =\p_t v_0(x,T) = 0$  and $v_k(x,0) =\p_t v_k(x,0) = 0$ for $k=1,2$ and $x\in [a,b]$. %so that the products of the form $v_0v_1v_2$ become dense in $L^1(\Omega\times [0,T])$. This recovers $q$.
 
To recover information about $q$ from \eqref{eq:integral_density}, we need to make an appropriate choice of the functions $v_1$ and $v_2$. Let $\chi\in C_c^\infty(\R)$ be a cut-off function supported close to $0\in \R$ with $\chi(0)=1$. For $\tau >0$, we let $H^\tau\in C_c^{\infty}(\mathbb{R})$ be a function defined by 
\[
H^\tau(\s l\s)=\chi(\s l\s)\tau^{1/2} \e^{-\frac{1}{2}\tau\, l^2}, \quad l\in \R.
\] 
% Let \[
% H(l)=\chi_\alpha(l)\tau^{1/2} \e^{-\frac{1}{2}\tau\, l^2}, \quad l\in \R.
% \]
Let $(x_0, t_0)\in [a,b]\times [0,T]$ and define a family of functions by
\begin{equation}\label{eq:H_functions}
\begin{split}
H_1^{\tau}(x,t)&:=H^\tau\big((x-x_0)-(t-t_0)\big), \\
H_2^{\tau}(x,t)&:= H^\tau\big((x-x_0)+(t-t_0)\big) 
\end{split}
\end{equation}
parametrized by $\tau$. 
Note that the functions $H^{\tau}_k(x,t)=0$, $k=1,2$, solve the wave equation and they are supported outside neighborhoods of $\{t=0\}$ and $\{t=T\}$ if the support of $\chi$ is chosen small enough. %We now define the functions $v_0,v_1$ and $v_2$. 
For $\tau>0$, we then set
\begin{equation}\label{eq:choice_of_fs}
\begin{cases}
 v_1=H_1^{\tau},& \ f_1=\big(H_1^{\tau}(a,t),H_1^{\tau}(b,t)\big)\\%= H_1^{\tau}|_\Sigma \\
 v_2=H_2^{\tau},& \ f_2=\big(H_2^{\tau}(a,t),H_2^{\tau}(b,t)\big)% H_2^{\tau}|_\Sigma.
 \end{cases}
 \end{equation}
 Here we suppressed the $\tau$ dependence of $v_1$ and $v_2$ in the notation. 
 As the function $v_0$ we take $H_1^{\tau}$ at fixed $\tau=1$:
 \begin{equation}\label{eq:v_0}
 v_0=H_1^{\tau}|_{\tau=1}, \quad \psi=v_0|_\Sigma.
 \end{equation}

\begin{figure}[ht!]
\centering
\includegraphics[scale=1.3]{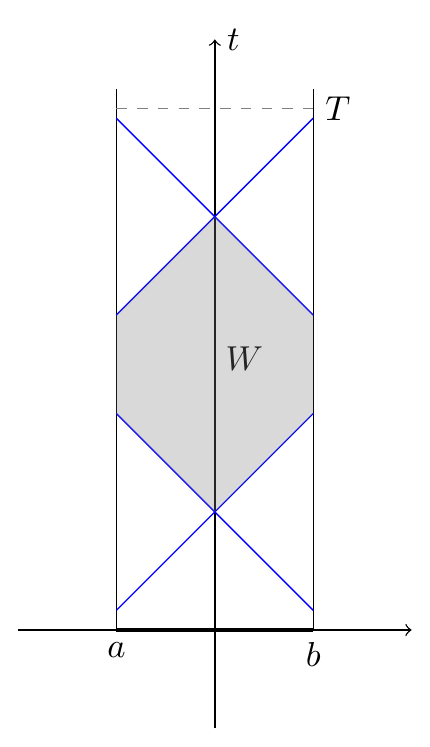}
\caption{ The set $W$, where we can recover the potential $q(x,t)$. In the 1+1 dimensional case this is the set, where we can make two approximate plane waves (one travelling to the left and the other to the right) collide. \label{fig:admissibleset}}
\end{figure}

Next we observe that the product of $v_1=H_1^\tau$ and $v_2=H_2^\tau$ approximates the delta distribution at $(x_0,t_0)$, when $\tau$ is large. Indeed, by~\cite[Lemma 3]{LLPT20} we have that
\begin{equation}\label{eq:tau2}
 q(x_0,t_0)-\frac{1}{\pi}\int_{[a,b]\times [0,T]} q(x,t) H_1^{\tau}(x,t)H_2^{\tau}(x,t)dx dt =\mathcal{O}(\tau^{-1/2}). %\longrightarrow \delta_{(x_0,t_0)}(x,t) \text{ as } \tau\to \infty,
\end{equation}
It was also proved in~\cite{LLPT20} that the implicit constant on the right hand side of \eqref{eq:tau2} is independent of the chosen point $(x_0, t_0)$. In particular $\frac{1}{\pi}H_1^{\tau}H_2^{\tau}$ converges to the delta  distribution at $(x_0,t_0)$ in the sense of distributions, as $\tau\to \infty$.
%Here $\delta_{(x_0,t_0)}(x,t)$ is the delta distribution at the point $(x_0,t_0)$. 
 %The required initial conditions for $v_0$, $v_1$ and $v_2$ are satisfied if $\alpha$ is chosen small enough. 
%  
%  as $v_0$ we take $H_1^{\tau,(x_0,t_0)}|_{\tau=1}$.
% \]
% v_k=H_k^{\tau,(x_0,t_0)}$, $k=1,2$, and for $v_0$ we take $H_1^{\tau,(x_0,t_0)}|_{\tau=1}$. 
It follows by taking $\tau\to \infty$ in~\eqref{eq:integral_identity_derivitve} that
\begin{equation}\label{eq:recovery_eq}
 -2\pi q(x_0,t_0)=\lim_{\tau\to\infty}\int_{\Sigma}\psi\s \p^2_{\eps_1\eps_2}|_{\eps_1=\eps_2=0}\Lambda_q(\eps_1f_1+\eps_2f_2)\s \d \Sigma.
\end{equation}
Since the right hand side of this equation is given by the DN map of \eqref{eq:intro_wave-eq}, we have recovered $q(x_0,t_0)$. Since $(x_0,t_0)$ in the above argument was arbitrary, we may similarly recover $q$ at any point in the set pictured in Figure \ref{fig:admissibleset}. 
This set is called the \emph{admissible set} and denoted by $W$. It is a set that can be reached from both the left and right sides of the lateral boundary by sending waves and from which waves can propagate back to the lateral boundary $\Sigma$.

\subsection{The numerical recovery method}\label{sec:num_rec}
Let us then explain how the theoretical recovery method described above transforms into a numerical recovery method. We conclude this section by explaining how we tested the quality of the numerical recovery method. 

Our numerical recovery method is based on the 
%We will use the 
identity \eqref{eq:recovery_eq}, which entails choosing suitable boundary values, taking a limit in $\tau$ and calculating derivatives with respect to the $\eps_1$ and $\eps_2$ parameters.   
% \noindent\textbf{1.} Fix a point $(x_0,t_0)$ where we aim to recover $a(x_0,t_0)$. Fix also the values of $\tau$ and $\eps$.
%  
% % \noindent\textbf{2.} Solve the equation \eqref{eq:intro_wave-eq} numerically with three different boundary values
% % \begin{align}
% %  \eps_1f_1+\eps_2f_2, \quad \eps_1f_1 \ \text{ and } \ \eps_2f_2.
% % \end{align}
% %%to obtain the normal derivatives of $u_{\eps_1f_1+\eps_2f_2}$, $u_{\eps_1f_1}$ and $u_{\eps_2f_2}$ on the lateral boundary $\Sigma=\{a,b\}\times[0,T]$.  
% %Here $f_1$ and $f_2$ are as in \eqref{eq:choice_of_fs}. This is done in Section \ref{}
% 
% \noindent\textbf{2.} Estimate the derivative
% \[
%  \p^2_{\epsilon_1, \epsilon_2}\big|_{\eps_1=\eps_2=0} \Lambda(\eps_1f_1+\eps_2f_2)
% \]
% for example by using finite differences \eqref{eq:data_finite_difference}.
%  
%  \noindent\textbf{3.} Integrate the finite differences in the time variable $t$ over $[0,T]$ to obtain an approximation of $a(x_0,t_0)$ by the formula \eqref{est:2D_noise}.
% %\end{enumerate}
%as the basis for the numerical recovery of the potential $q$ from the DN map. 
To evaluate the mixed derivative $\p^2_{\eps_1\eps_2}|_{\eps_1=\eps_2=0} \Lambda_q(\eps_1f_1+\eps_2f_2)$ in \eqref{eq:recovery_eq}, we experimented with two approaches to evaluate this derivative: 

%\noindent
\textbf{(1)} A finite difference approximation. 

%\noindent
\textbf{(2)} An approach using regularization, where differentiation is regarded as an inverse operation to integration.

\noindent Let us first discuss the first approach \textbf{(1)}. Finite differences were used in \cite{LLPT20} to prove a H\"older stability estimate for the recovery of $q$. There the mixed derivatives $\p^2_{\eps_1\eps_2}|_{\eps_1=\eps_2=0}\Lambda_q(\eps_1f_1+\eps_2f_2)$
%$\p^2_{\eps_1\eps_2}|_{\eps_1=\eps_2=0}\s u_{\eps_1f_1+\eps_2f_2}$ 
were replaced by the finite differences 
\begin{equation}\label{eq:data_finite_difference}
D^2_{\epsilon_1, \epsilon_2}\big|_{\eps_1=\eps_2=0} \Lambda_q(\eps_1f_1+\eps_2f_2):= \frac{1}{\epsilon_1 \epsilon_2} \left(\Lambda_q(\eps_1f_1+\eps_2f_2)- \Lambda_q(\eps_1f_1)-\Lambda_q(\eps_2f_2)\right).
\end{equation}
%Here $u_{\eps_1f_1+\eps_2f_2}$ solves~\eqref{eq:prelim_wave-eq_recovery} as before. 
For finite differences, the integral identity % instead of mixed derivative, to have that $w$ satisifies
%  \[
%  w=D_{\eps_1\eps_2}^2\Big|_{\eps_1=\eps_2=0} u_{\eps_1 f_1 + \eps_2 f_2}+ \mathcal{O}(\langle\eps_1,\eps_2\rangle^3).
% \]
\begin{equation}\label{eq:finite_difference_recovery}
\begin{aligned}
  -2\int_{[a,b]\times [0,T]} q\s v_0\s v_1\s v_2\s \d x \s \d t &=\int_{\Sigma}v_0\s D_{\eps_1,\eps_2}^2\Big|_{\eps_1=\eps_2=0}\Lambda_q(\eps_1f_1+\eps_2f_2) \s \d S\\
  &\qquad \qquad\qquad\qquad + \frac{1}{\eps_1\eps_2}\int_{[a,b]\times [0,T]} v_0\square\s \tildeR\s \d x \s \d t
  \end{aligned}
 \end{equation}
 is the replacement of \eqref{eq:integral_identity_derivitve}. This identity holds for any $f_1,f_2\in H^s(\Sigma)$, $s\in \N$, and corresponding solutions $v_1$ and $v_2$ to \eqref{eq:vj_equations}.
Here $\square \s \tildeR=\mathcal{O}(\langle\eps_1,\eps_2\rangle^3)$ in $H^{s+1}([a,b]\times [0,T])$, 
%where the same calculations were used to derive a stability result for the recovery of $q$. ,
 %for details see \eqref{eq:energy_norm} and \eqref{eq:tildeR}--\eqref{est:square_tildeR}. 
where $\langle\eps_1,\eps_2\rangle^3$ denotes an unspecified homogeneous polynomial of order $3$ in $\eps_1$ and $\eps_2$. We refer to~\cite{LLPT20} for more details. Here we suffice to remark that the last integral in \eqref{eq:finite_difference_recovery} is of the size $\max(\eps_1,\eps_2)$, which is small. 
%For the choices of $v_0$, $v_1$ and $v_2$ given in previous section, the left hand side of \eqref{eq:finite_difference_recovery} is approximately the value of $q$ at $(x_0,t_0)$ (times $-2\pi$).

%Recall from \eqref{eq:tau2} that for our choices of $v_0$, $v_1$ and $v_2$, the left hand side of \eqref{eq:finite_difference_recovery} is approximately the value of $q$ at $(x_0,t_0)$ times $-2\pi$.

We choose $f_1,f_2$ and $\psi$ as in \eqref{eq:choice_of_fs}--\eqref{eq:v_0}. We also introduce noise $\mathcal{E}$ to our measurements. We assume that $\mathcal{E}$ is a bounded, 
\begin{equation}%\label{eq:size_of_noise}
  \norm{\E (f)}_{H^r(\Sigma)}\leq \delta, 
 \end{equation}
possibly non-linear, mapping $H^{s+1}(\Sigma) \to H^r(\Sigma)$, $r\in \R$ and $r \leq s$.
In this case, for our choices for $f_1$ and $f_2$, the equation~\eqref{eq:finite_difference_recovery}  is replaced by
\begin{equation}\label{eq:finite_difference_recovery_with_noise}
\begin{aligned}
  -2\int_{[a,b]\times [0,T]} q\s v_0\s v_1\s v_2\s \d x \s \d t &=\int_{\Sigma}\psi\s D_{\eps_1, \eps_2}^2\Big|_{\eps_1=\eps_2=0}(\Lambda_q+\mathcal{E})(\eps_1f_1+\eps_2f_2) \s \d \Sigma\\
  &\qquad  + \frac{1}{\eps_1\eps_2}\int_{[a,b]\times [0,T]} v_0\square\s \tildeR\s \d x \s \d t.
  \end{aligned}
 \end{equation}
 %for any $f_1,f_2\in H^s(\Sigma)$.
%It was shown in~\cite[Theorem]{LLPT20} that by assuming \eqref{eq:size_of_noise} and from the identity \eqref{eq:finite_difference_recovery_with_noise} it follows that 
It was shown in \cite[Theorem 3]{LLPT20} that there is a constant $C>0$ such that
 \begin{equation}\label{est:2D_noise}
\begin{aligned}
\Big|- q(x_0, t_0) -\frac{1}{2\pi}\int_{\Sigma} \psi\,D_{\eps_1,\eps_2}^2\Big|_{\eps_1=\eps_2=0}(\Lambda_q + \E)(\eps_1f_1 +\eps_2f_2) \d\Sigma 
\Big|   \leq C \delta^{\sigma (s)},
\end{aligned}
\end{equation}
The exponent $\sigma(s)$ is $\frac{1}{6s+15}$. The estimate~\eqref{est:2D_noise} is a H\"older stability estimate for the recovery of $q$. It is obtained by optimizing the right hand side of~\eqref{eq:finite_difference_recovery_with_noise} with respect to the parameters $\eps_1,\eps_2$ and $\tau$. If one computes the integral in \eqref{est:2D_noise} and if the noise is small, then the integral is approximatively the potential $q$ at the point $(x_0,t_0)$. This is the basis of our numerical reconstruction of $q$.

The optimal choice of $\eps_1,\eps_2$ and $\tau$ depends on some constants, mainly on a priori assumptions on the potentials and the size of the noise. While the constants can in principle be calculated, we do not attempt to do that. Instead, we suffice to experiment by choosing small non-zero values for $\eps_1,\eps_2>0$ and a large  finite value for $\tau$. 
%In this paper we set $\eps:=\eps_1=\eps_2$, and just choose $\eps>0$ and $\tau>0$ so that $\eps$ is very small compared to $\tau$ and so that we get good numerical reconstruction. 
We will discuss our choices for $\eps_1$ and $\eps_2$ more carefully in Sections~\ref{sec:Dregularization} and~\ref{sec:tau_and_eps}. The large parameter $\tau$ will be chosen by studying the size of the error term in \eqref{eq:tau2}.
 We observe experimentally that $\tau=700$ and $\max(\eps_1,\eps_2)\approx 0.01$ produces good results in the examples we consider in Section \ref{sec:examples}.

%The paper \cite{LLPT20} also optimizes \eqref{eq:finite_difference_recovery} with respect to $\eps_1$ and $\eps_2$ and  the parameter $\tau$ implicit in the solutions $v_1$ and $v_2$. However the optimization contains some constants, which are hard or laborsome to evaluate in practice. In this paper, we suffice to experiment by choosing small \f{In what sense small?} non-zero values for $\eps_1,\eps_2>0$ and a large \f{In what sense large?} finite value for $\tau$. The equation \eqref{eq:finite_difference_recovery}  then yields the value of $q$ approximatively at $(x_0,t_0)$.

Let us then describe our second approach \textbf{(2)} for calculating the mixed derivative $\p^2_{\eps_1\eps_2}|_{\eps_1=\eps_2=0} \Lambda_q(\eps_1f_1+\eps_2f_2)$. For this, let us denote
\[
 F(s,t)=\Lambda_q(s f_1+tf_2).
\]
By the fundamental theorem of calculus
\begin{align}\label{eq:fund_thm2_intro}
\begin{split}
\int_\beta^{\eps_1}\int_\alpha^{\eps_2} \p_s\p_tF(s,t) \d s\s  \d t
&= F(\eps_1,\eps_2) - F(\alpha,\eps_2) - F(\eps_1,\beta) + F(\alpha,\beta) \\
& \quad \text{ for all } \eps_1 \text{ and } \eps_2 \text{ small enough}.
\end{split}
\end{align}
Since the identity above holds for all $\eps_1$ and $\eps_2$ small enough, we consider it as system of equations parametrized by $\eps_1$ and $\eps_2$. Assume now that there is some error or noise on the right hand side of the equations in \eqref{eq:fund_thm2_intro}. Fix also $\alpha$ and $\beta$. The approach \textbf{(2)} is based on the question if in this case we can find an approximation for $\p_s\p_tF(s,t)$ by solving the system \eqref{eq:fund_thm2_intro}? The system is at least typically ill-posed, so in general the answer is no. However, one can choose a suitable regularization scheme depending on known properties of the measurement noise and the boundary data to solve the system approximately. We find an approximate solution $\p_s\p_tF(s,t)$ by using Tikhonov regularization. An approximation of $\p^2_{\eps_1\eps_2}|_{\eps_1=\eps_2=0}\Lambda(\eps_1f_1+\eps_2f_2)$ is then $F(0,0)$. This is the approach \textbf{(2)}. We call the approximate solution obtained in this way the \emph{regularized mixed derivative}.
% 
% we will show that an approximation to the second derivative $\p^2_{\eps_1\eps_2}|_{\eps_1=\eps_2=0}\Lambda(\eps_1f_1+\eps_2f_2)$ can be found by approximately solving for a function $F$ the family of equations
% \begin{equation}%\label{eq:fund_thm2}
% \begin{split}
% \int_c^{\eps_1}\int_a^{\eps_2} F(s,t) ds dt
% &= \Lambda(\eps_1,\eps_2) - \Lambda(a,\eps_2) - \Lambda(\eps_1,c) + \Lambda(a,c)
% \end{split}
% \end{equation}
% parametrized by $\eps_1$ and $\eps_2$ small. Of course, if $F(s,t)=\p^2_{st}\Lambda(sf_1+tf_2)$ then the above identity holds. The point here is that when there is error or noise on the right hand side of the identity, we can still to find an approximation for $\p^2_{st}\Lambda(sf_1+tf_2)$
% 
% An approximation of $\p^2_{\eps_1\eps_2}|_{\eps_1=\eps_2=0}\Lambda(\eps_1f_1+\eps_2f_2)$ is then $F(0,0)$. If 
% 
% . Indeed, one can think of finding a derivative as the inverse problem to calculating an integral. This approach has the benefit that one can choose a suitable regularization scheme depending on known properties of the measurement noise and the boundary data.\footnote{Rewrite this more clearly?}

We end this section by discussing how our reconstruction method is implemented. The reconstruction is non-iterative and splits into a few steps. In Section \ref{sec:examples} we compare exact potentials $q$ to the approximative ones obtained by our reconstruction method. 
%Numerical examples are given in Section \ref{sec:examples}.
  To demonstrate reconstruction of a potential function, we do the following steps:
%Given boundary values $\Lambda(\eps_1f_1+\eps_2f_2)$, the steps are:\footnote{Rewrite this.}
\\
%\begin{enumerate}e 
 %\item 
 \noindent\textbf{a)} Fix a grid for the domain $[a,b]\times [0,T]$ and a point $(x_0,t_0)$ of the grid. We aim to recover $q(x_0,t_0)$. Corresponding to $(x_0,t_0)$, let us choose the boundary values $f_1$ and $f_2$ as in \eqref{eq:choice_of_fs}. Fix also a discrete set for the small values for $\eps_1$ and $\eps_2$ and fix also $\tau$ large.  \\
 \noindent\textbf{b)} Solve the forward problem \eqref{eq:intro_wave-eq} by using the boundary values $\eps_1f_1+\eps_2f_2$ by using the finite difference approach explained in Section \ref{sec:forward}. Compute the boundary normal derivatives of the corresponding solutions $u_{\eps_1f_1+\eps_2f_2}$ on the lateral boundary using first order central differences. 
 
 %Fix also the values of $\tau$ and $\eps$.
 
% \noindent\textbf{2.} Solve the equation \eqref{eq:intro_wave-eq} numerically with three different boundary values
% \begin{align}
%  \eps_1f_1+\eps_2f_2, \quad \eps_1f_1 \ \text{ and } \ \eps_2f_2.
% \end{align}
%%to obtain the normal derivatives of $u_{\eps_1f_1+\eps_2f_2}$, $u_{\eps_1f_1}$ and $u_{\eps_2f_2}$ on the lateral boundary $\Sigma=\{a,b\}\times[0,T]$.  
%Here $f_1$ and $f_2$ are as in \eqref{eq:choice_of_fs}. This is done in Section \ref{}

\noindent\textbf{c)} Estimate the mixed derivative numerically
\[
 \p^2_{\epsilon_1\s \epsilon_2}\big|_{\eps_1=\eps_2=0} \s \Lambda_q(\eps_1f_1+\eps_2f_2)
\]
by using either the finite difference $D_{\eps_1,\eps_2}^2|_{\eps_1=\eps_2=0}(\Lambda_q + \E)(\eps_1f_1 +\eps_2f_2)$ in the approach \textbf{(1)} or the regularized mixed derivative in the approach \textbf{(2)}. Numerical implementation of the latter  is given in Section \ref{sec:Dregularization}.
 
 \noindent\textbf{d)} Compute the integral in \eqref{est:2D_noise}, where $D_{\eps_1,\eps_2}^2\Big|_{\eps_1=\eps_2=0}(\Lambda_q + \E)(\eps_1f_1 +\eps_2f_2)$ is replaced by the regularized derivative in the case \textbf{(2)}. This results in an approximation of $q(x_0,t_0)$. This part is done in Section \ref{sec:IP}. Repeat the now described steps at all grid points inside the admissible set to obtain an approximation of $q$. (The admissible set was illustrated in Figure \ref{fig:admissibleset}.)

\section{Forward model}\label{sec:forward}

We use a finite difference approach to numerically solve the nonlinear wave equation \eqref{eq:intro_wave-eq}. For a discussion about the finite difference method for linear wave equations, see e.g. \cite{MG,Smith}.
Let $N_x$ and $N_t$ be positive integers.
We discretize the problem \eqref{eq:intro_wave-eq} by dividing the spatial and time-domains $[a,b]$ and $[0,T]$ into uniformly spaced grids
\begin{align*}
x_i &= a + (i-1)\Delta x, \quad i=1,\ldots, N_x+1,\\
t_j &= (j-1)\Delta t, \quad j=1,\ldots, N_t+1,\\
\end{align*}
of length $\Delta x = (b-a)/N_x$ and $\Delta t = T/N_t$.
As we assume that the wave speed is constant $1$ within the domain, we have the Courant-Friedrichs-Lewy (CFL) number $c = \Delta t /\Delta x$.
We always choose the numbers $N_x$ and $N_t$ so that $c<1$.
The equation \eqref{eq:intro_wave-eq} is represented in the discretized form as
\begin{equation}\label{eq:FDTD}
\begin{split}
&\frac{u(x_i,t_{j+1}) - 2u(x_i,t_j) + u(x_i,t_{j-1})}{\Delta t^2} \\
&\qquad\qquad\qquad\qquad=\frac{u(x_{i+1},t_j) - 2u(x_i,t_j) + u(x_{i-1},t_j)}{\Delta x^2} + a(x_i,t_j)u(x_i,t_j)^2.
\end{split}
\end{equation}
%
%The equation \eqref{eq:FDTD} converges to \eqref{eq:intro_wave-eq} as $\Delta x,\Delta t\to 0$. Thus the scheme is consistent. It is also accurate to second order with respect to the truncation error. This means that a smooth solution to \eqref{eq:intro_wave-eq} satisfies \eqref{eq:FDTD} up to an error of size $\Delta t^2+\Delta x^2$.
A smooth (for example $C^4$) solution to \eqref{eq:intro_wave-eq} satisfies \eqref{eq:FDTD} up to an error of size $\leq C(\Delta t^2+\Delta x^2)$, which means that the  finite difference scheme \eqref{eq:FDTD} is consistent with \eqref{eq:intro_wave-eq} and accurate to second order.
To model vanishing Cauchy data, we initialize
$$
u(x_i,t_1) = u(x_i,t_2) = 0
$$
for all $x_i$, $i=1,\ldots, N_x+1$. %This corresponds to vanishing Cauchy data. 
Similarly, the boundary values are enforced by setting
$$
u(a,t_j) = f^L(t_j) \text{ and } u(b,t_j) = f^R(t_j),
$$
for all $t_j$, $j=1,\ldots,N_t+1$.

To evaluate the Dirichlet-to-Neumann map on the lateral boundary $\Sigma=\{a,b\}\times [0,T]$ we use first order central differences for the normal derivatives as
\begin{equation}\label{eq:Neumann_derivatives}
\begin{split}
\partial_\nu u(a,t_j) &\approx -\frac{u(x_3,t_j) - u(x_1,t_j)}{2\Delta x} \\% =: -D u(x_1,t_j),\\
\partial_\nu u(b,t_j) &\approx \frac{u(x_{N_x+1},t_j) - u(x_{N_x-1},t_j)}{2\Delta x}% =: D u(b,t_j)
\end{split}
\end{equation}
for all $t_j$, $j=1,\ldots, N_t+1$.

\subsection{Example}\label{sec:forward_example}
To verify the convergence of the forward solver experimentally we calculated examples of numerical solutions $u$ of \eqref{eq:FDTD} in the domain $[a,b]\times[0,T]=[-0.5,0.5]\times[0,3]$ using increasingly denser computation grids. 

\begin{figure}[ht!]
\centering
\includegraphics[width = 0.8\textwidth]{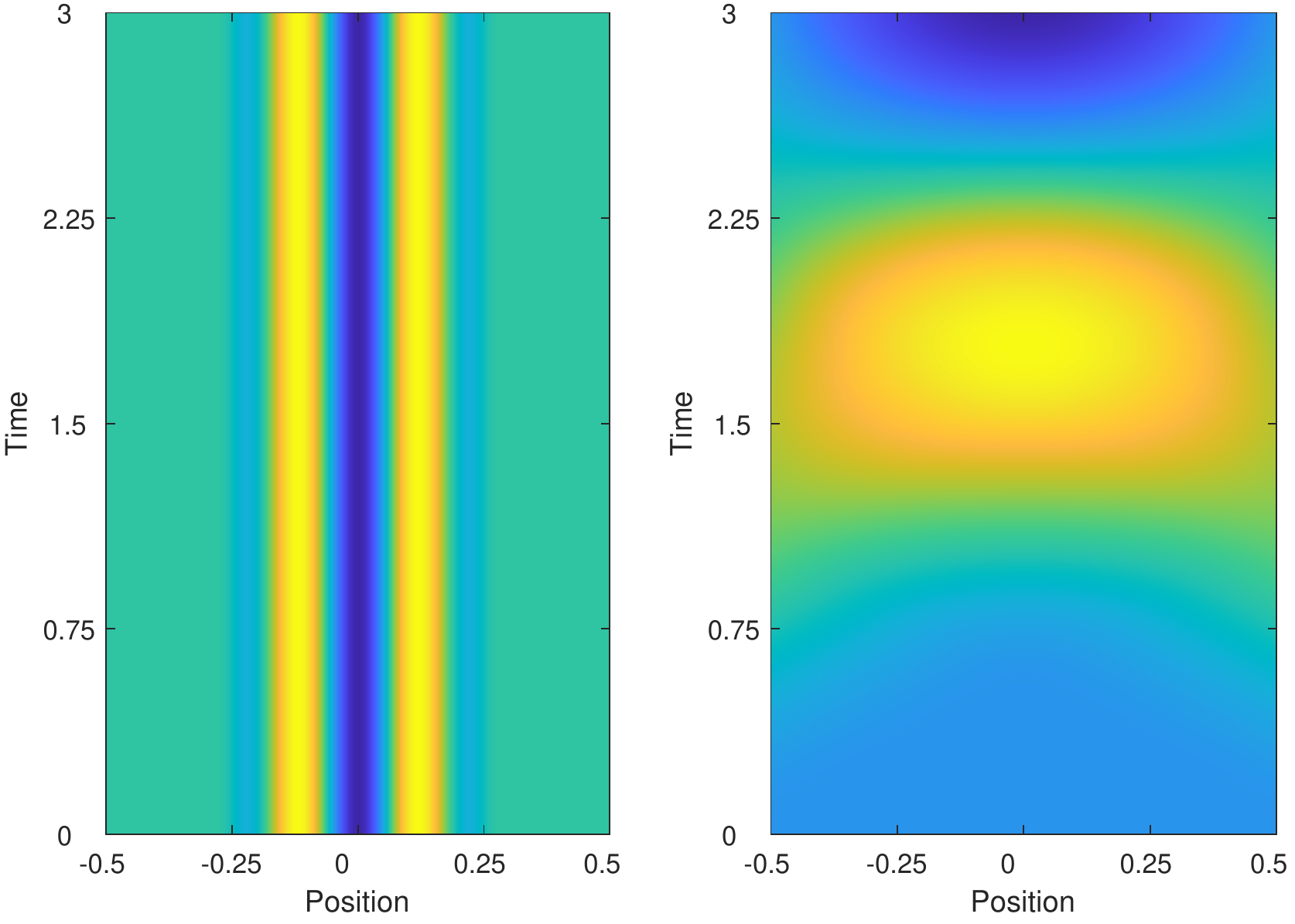}
\caption{The example of Section~\ref{sec:forward_example}. The potential function $q(x,t)$  (left) and the corresponding solution $u$ of \eqref{eq:intro_wave-eq} (right).  \label{fig:u_sinsquared}}
\end{figure}

Let us describe one example here. Let $\varphi\in C_c^\infty(\R)$ be the bump function
\begin{equation}\label{eq:testfun}
\varphi(s) = 
\begin{cases}
\exp\left(\frac{1}{s^2-1}\right), & \text{when } s^2<1,\\
0, & \text{otherwise}.
\end{cases}
\end{equation}
Let the potential function $q$ be
\begin{equation}\label{eq:example_a}
q(x,t) = \varphi(10x/3)\sin(x)
\end{equation}
and let the boundary values of \eqref{eq:intro_wave-eq} be
\begin{equation*}
f^L(t) = f^R(t) = 0.01 \sin(t)^2.
\end{equation*}
The numerical solutions $u_h$ of \eqref{eq:FDTD}, with $q$ and $f^L$ and $f^R$ as above, were evaluated with $\Delta x = 1/(4\times 2^h)$ and $\Delta t = 1/(32\times 2^h)$ at the grid points $x_i = -0.5 + (i-1)\Delta x$ and $t_j = (j-1)\Delta t$, where $i=1,\ldots,4\times 2^h+1$, $j=1,\ldots,32\times 2^h+1$, and $h=2,\ldots,10$. The numerical solution $u_{h=10}$ is depicted in Figure~\ref{fig:u_sinsquared}.

\begin{figure}[ht!]
\centering
\includegraphics[width = 0.9\textwidth]{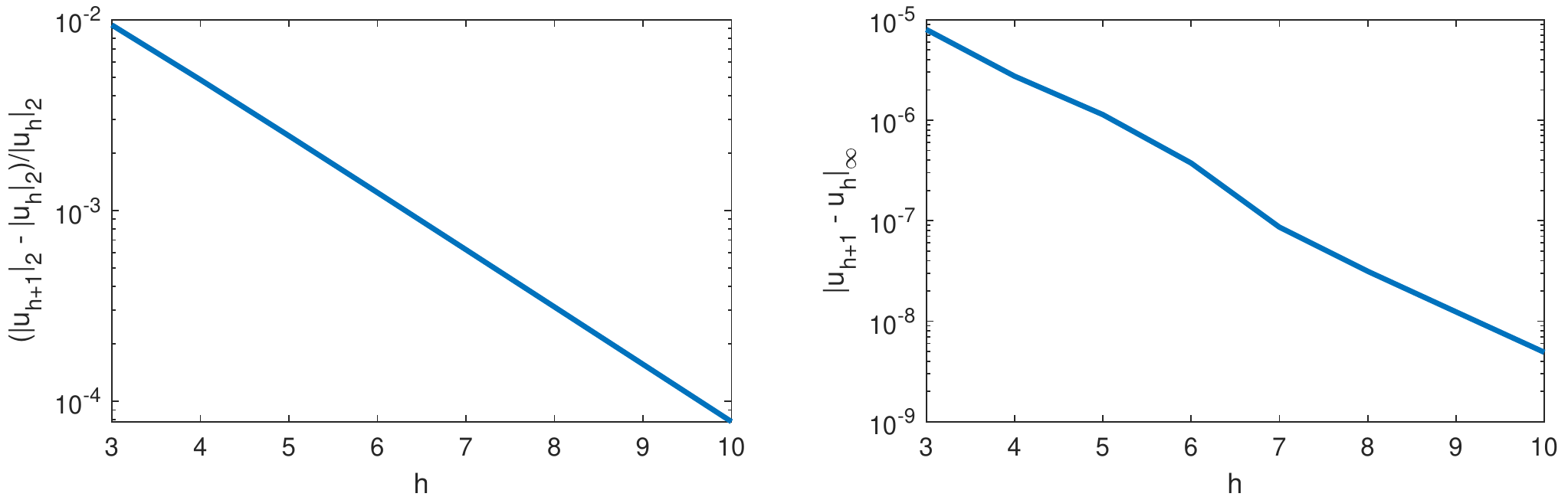}
\caption{Left: The relative difference $(| u_{h+1} |_2 - | u_{h} |_2)/| u_{h} |_2$ between $L^2$-energies of solutions $u_h$ to \eqref{eq:FDTD} at successively finer grids in log-scale. Right: the absolute differences $\Vert u_{h+1}-u_h\Vert_{L^\infty}$ evaluated at the common grid points between the successively finer grids in log-scale. We see that the $L^2$-energy stays bounded along with finer grids and the values of the numerical solutions converge at the common grid points.
 \label{fig:convergence}}
\end{figure}

To study the convergence of the numerical scheme, we recorded the $L^2$-energy
$$
| u_h |_{2} := \sqrt{\sum_{i,j} u_h(x_i,t_j)^2\Delta x\Delta t}
$$
of each solution $u_h$ and
%
%We also evaluated the absolute differences \f{Remove the below.}
%$$
%| u_{h+1}(x_i,t_j) - u_h(x_i,t_j)|
%$$
%at the common grid points between the two successive iterations $u_{h+1}$ and $u_h$. 
%%
the maximum of the absolute differences 
$$
| u_{h+1} - u_h |_{\infty} := \max_{i,j} | u_{h+1}(x_i,t_j) - u_h(x_i,t_j)|
$$
at the common grid points.
The maximum absolute differences and the relative differences
$$
\frac{| u_{h+1} |_2 - | u_{h} |_2}{| u_{h} |_2}
$$
between the $L^2$-energies of successive solutions are show in Figure~\ref{fig:convergence}.

\section{Regularized mixed derivatives}\label{sec:Dregularization}
As explained in Section \ref{sec:num_rec}, to recover the potential $q(x,t)$ numerically from the noisy DN map $\widetilde\Lambda$,
$$
\widetilde\Lambda(\eps_1 f_1+ \eps_2 f_2) := \Lambda_q(\eps_1 f_1+ \eps_2 f_2) + \mathcal{E},
$$
we need to evaluate the mixed derivative 
$$
 \p^2_{\epsilon_1\s \epsilon_2}\big|_{\eps_1=\eps_2=0} \s\widetilde\Lambda(\eps_1 f_1+ \eps_2 f_2).
$$
%Here we use the shorthand notation $\widetilde\Lambda$ for a  noisy measurement:
% 
As is well known (and easy to see), numerical differentiation is highly sensitive to noise. Before proceeding to our actual inverse problem, in this section we consider the approach \textbf{(2)} explained in Section \ref{sec:num_rec} to compute mixed derivatives approximatively using regularization. We call the result \emph{regularized mixed derivative}.

%authors have noted that the problem of evaluating numerical derivatives is an inverse problem to calculating integrals.~\f{Let's be more specific.} 
%
%As is often the case with inverse problems, evaluating $g = f'$ is an ill-posed problem, since $g$ does not depend continuously on $f$.~\f{What is this?}
%
Our method derives from the results by Cullum \cite{Cu1971}, who used Tikhonov regularization to estimate numerical derivatives.
For other recent results on numerical differentiation of noisy data we refer to the papers \cite{Ch2011,KR2014} and the references therein.

Our regularized method to compute $\p^2_{\epsilon_1 \epsilon_2}\big|_{\eps_1=\eps_2=0} \widetilde\Lambda(\eps_1f_1+\eps_2f_2)$ is as follows.  For this, we consider the problem of finding the mixed derivative $\partial_s\partial_tg(s,t)$ of a generic function   
%\footnote{fix some nice notation. change $[a,b]$.} 
% 
% As explained in Section \ref{sec:num_rec}, to recover the potential $q(x,t)$ numerically from the noisy DN map $\widetilde\Lambda$,
% $$
% \widetilde\Lambda(\eps_1 f_1, \eps_2 f_2) := \Lambda(\eps_1 f_1, \eps_2 f_2) + \mathcal{E},
% $$
% we need to evaluate the mixed derivative ~\f{Finite differnce of actual derivative here?}
% $$
% D^2_{\eps_1=\eps_2=0} \widetilde\Lambda(\eps_1 f_1, \eps_2 f_2).
% $$
% %Here we use the shorthand notation $\widetilde\Lambda$ for a  noisy measurement:
% %
% Let us consider the problem of evaluating the two-dimensional derivative $D^2_{xy} = \partial_1\partial_2$ numerically.\footnote{fix some nice notation} 
%
%Let 
$g\in C^2([\alpha,\alpha_0]\times [\beta,\beta_0])$. 
By the fundamental theorem of calculus
\begin{equation}\label{eq:fund_thm1}
\int_\alpha^x\p_1 g(s,y) \d s = g(x,y)-g(\alpha,y),\quad x\in [\alpha,\alpha_0].
\end{equation}
For $(x,y)\in [\alpha,\alpha_0]\times [\beta,\beta_0]$, it then follows that
\begin{equation}\label{eq:fund_thm2}
\begin{split}
\int_\beta^y\int_\alpha^x \p_1\p_2 g(s,t) \d s \s \d t
&=
\int_\beta^y \partial_2 \big(
g(x,t) - g(\alpha,t)
\big)\d t
\\
&= g(x,y) - g(\alpha,y) - g(x,\beta) + g(\alpha, \beta),
\end{split}
\end{equation}
%for $(x,y)\in [\alpha,\alpha_0]\times [\beta,\beta_0]$. 

% Assume now that the discrete measurements~\f{Range of the indices i,j? Remove below line. }
% $$
% f_{i,j}:= f(x_i,y_j) + \eps_{i,j},
% $$
% $i = 1,2\ldots, N_i\in \N$ and $j= 1,2\ldots, N_j\in \N$
% are corrupted by an identically distributed~\f{Is any noise ok at this point? I.d. wrt to i,j} noise $\eps_{i,j}$ with zero expectation value.
%

 Let $N_i$ and $N_j$ be positive integers. Let us consider a grid obtained by choosing $N_i$ and $N_j$ points from $[\alpha,\alpha_0]$ and $[\beta,\beta_0]$ respectively. 
%
%We discretize the problem \eqref{eq:intro_wave-eq} by dividing the spatial and time-domains $[a,b]$ and $[0,T]$ into uniformly spaced grids
Let us denote the discrete values of $\p_s\p_t f(s,t)$ at the grid points by a vector $\y$. The vector $\y$ has thus length $N_iN_j$. Let us also simply denote by $g$ the vector obtained from the values of $g$ at the grid points. 
%
%Adding noise to the the values of $f(s,t)$ at the grid points, 
By considering integration as a Riemann sum, the equation \eqref{eq:fund_thm2} in discretized form reads
%data~\f{What do we mean by measurements here?} $f$ as
\begin{equation}\label{eq:discrete_derivative}
A\y =  g,
\end{equation}
where the matrix $A$, which is of the size $(N_iN_j)\times (N_iN_j)$, is called the anti-differentiation matrix. %It is a discretized form of the double integral in \eqref{eq:fund_thm2} considered as a Riemann sum. 
%Here also $f$ is understood as the vector obtained from $(f_{i,j})$ obtained from the values of $f$ at the grid points. 
%where $A$ is a matrix of the size $(N_iN_j)\times (N_iN_j)$. 

Let us then add noise $\eps_{i,j}$ to the data $g$ by setting 
$$
\g_{i,j} = g(x_i,y_j) - g(\alpha,y_j) - g(x_i,\beta) + g(\alpha,\beta) + \eps_{i,j}.
$$
In our inverse problem the noise is independent of the grid points and has zero expected value at each grid point. 
%Here the matrix $A$ is called \emph{anti-differentiation} matrix. It is a discretized form of the double integral in \eqref{eq:fund_thm2} considered as a Riemann sum.
%
The linear system 
\[
 A\widetilde\y = \g
\]
 can then be solved approximatively via a choice of a regularization method. The approximative solution will then be an approximation of $\p_s\p_tg(s,t)$ at the grid points. 
Since the solution to \eqref{eq:intro_wave-eq} is at least $C^1$, we can use this knowledge as an \emph{a priori} information to help us solve the linear system $A\widetilde\y = \g$.
Our choice is to use the generalized Tikhonov regularization method 
(see e.g. the books \cite{Hansen2,MS2012}) written as the minimization problem 
\begin{equation}\label{eq:tikhonov}
\y_\mathrm{reg} = \mathrm{arg\,min}_{\widetilde\y}
\big\{
\Vert A\widetilde\y - \g\Vert_{L^2} + \lambda \Vert \widetilde\y \Vert_{H^1}
\big\}, \quad \lambda>0.
\end{equation}
The regularized solution can be obtained by solving 
\begin{equation}\label{eq:linear_system}
(\lambda D^TD + A^TA)\y_\mathrm{reg} = A^T \g,\quad \lambda>0.
\end{equation}
Here the matrix $D$ determines the discretization of the Sobolev $H^1$-norm,
$$
D^T = [I, D_1^T, D_2^T],
$$
where $I$ is the identity matrix, and $D_1$ and $D_2$ are the difference matrices to the $x$- and $y$-directions respectively.
For the matrices $D_1$ and $D_2$ we used periodic differences.
The regularization parameter $\lambda>0$ is chosen experimentally.
(The equation \eqref{eq:linear_system} can be solved since $D^TD = I^TI + D_1^TD_1 + D_2^TD_2$ is positive definite and $A^TA$ is positive semi-definite.) 

\begin{figure}[ht!]
\centering
\includegraphics[width = 0.9\textwidth]{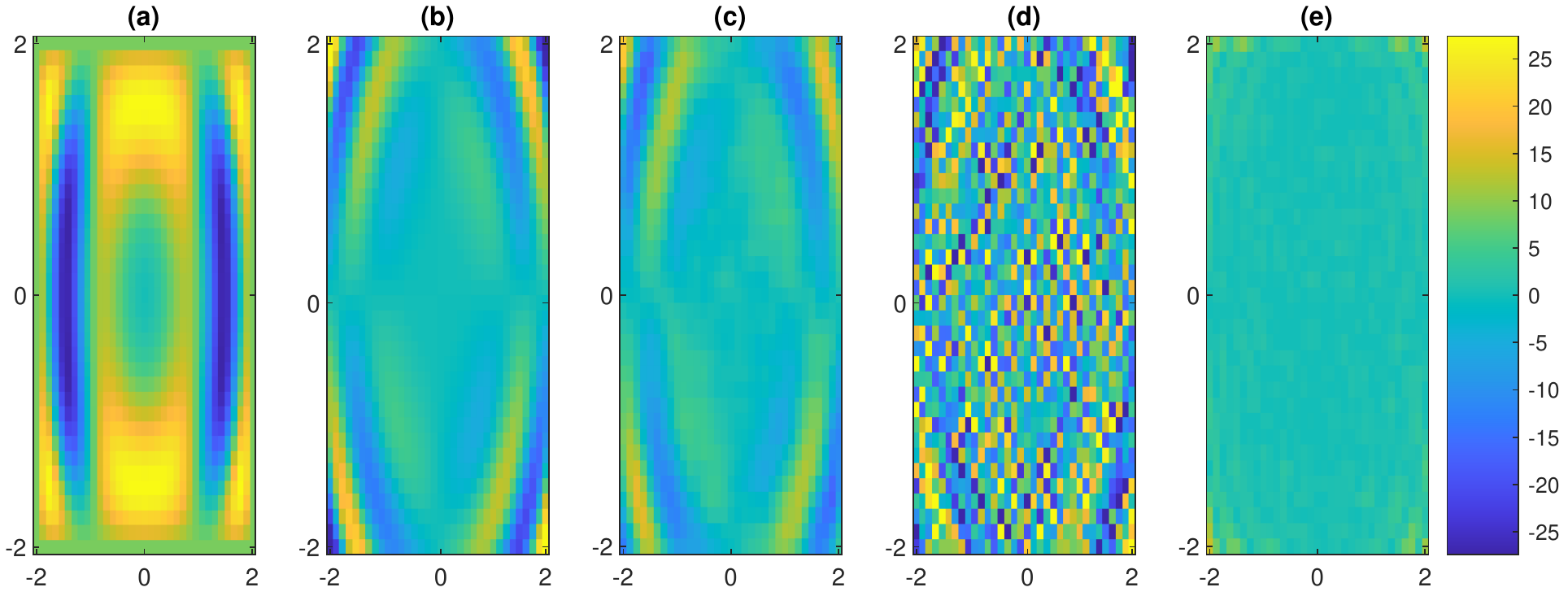}\\
\includegraphics[width = 0.9\textwidth]{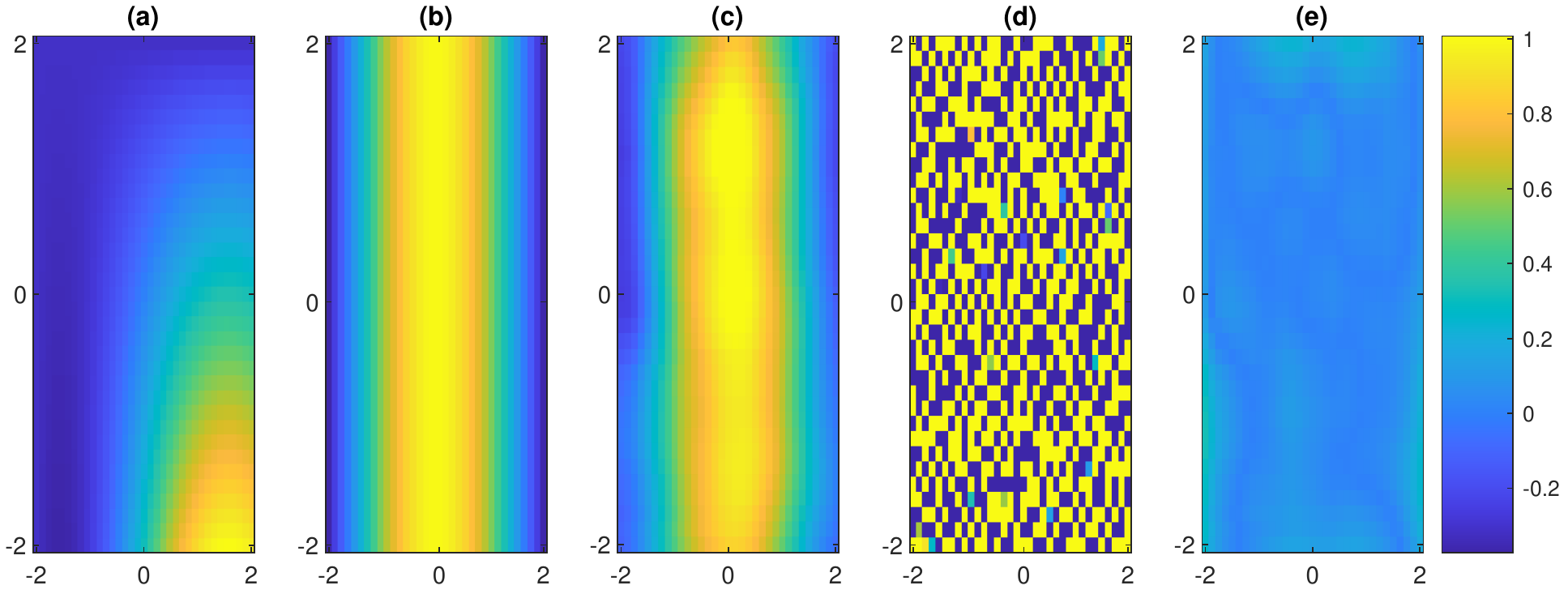}
\caption{Demonstration of regularization of numerical differentiation of noisy measurements. Example 1, the function \eqref{eq:ref_ex1}, on top row and example 2, the function in \eqref{eq:ref_ex2}, on bottom row. (a) The function $g$ to be differentiated. (b) The exact derivative $\p_x\p_yg$ of the function $f$. (c) Regularized approximation $\y_\mathrm{reg}$ to the derivative $\y$ of the noisy measurement $\g_{i,j}$. (d) Finite difference approximation of the derivative of the noisy measurement $g_{i,j}+\eps_{i,j}$. (e) The absolute difference $|\p_x\p_y g - \y_\mathrm{reg}|$. Experimentally chosen regularization parameters $\lambda=0.00001$ and $\lambda=0.01$ have been used in example 1 and 2 respectively. %, while for example 2 $\alpha=0.01$, both chosen experimentally.
\label{fig:regularizedFD}}
\end{figure}

To demonstrate this approach of numerical differentiation, we calculated regularized mixed derivatives of the functions
\begin{align}
g(x,y) &= \sin(x^2 + y^2) \label{eq:ref_ex1} \\
g(x,y) &= y\sin(x) \label{eq:ref_ex2}
\end{align}
over the region $(x,y)\in [-2,2]\times[-2,2]$.
The functions were sampled at $35^2 = 1225$ equally spaced grid points.
Both measurements were then corrupted by Gaussian noise with standard deviation $\sigma$ equal to $10\,\%$ of the maximum of the function $g$:
$$
\sigma=0.1\max_{(x,y)\in [-2,2]^2} |g(x,y)|.
$$
The resulting regularized derivatives $\p_x\p_y g(x,y)$ are depicted in Figure~\ref{fig:regularizedFD}.
To demonstrate the effects of regularization, we also included the finite difference approximations of the derivatives in the figures.
The effectiveness of regularization is visually striking.

\section{The inverse problem}\label{sec:IP}
Let $W$ be the admissible set defined in Section \ref{sec:prelim} and illustrated in Figure \ref{fig:admissibleset}.  This is a set which can be reached by sending waves from the lateral boundary $\Sigma$ and from which wave signals can also be detected on $\Sigma$. It is the diamond-shaped region in the $(x,t)$-plane given by the conditions
\begin{align*}
& t>x-a,\\ &t< -x+b+T,\\ &t>-x+b,\\ &t<x-a+T, \\
& x\in [a,b] \quad\text{and}\quad t\in [0,T].
\end{align*}
%Th is the set which we can reach by sending waves from the boundary $\Sigma$ and for which the corresponding measurements can also be detected on $\Sigma$.
%
Let $(x_0,t_0)\in W$. 

\smallskip
\noindent\textbf{Approach (1): Finite difference approximation of $\p_{\eps_1\eps_2}^2|_{\eps_1=\eps_2=0}\widetilde\Lambda$}.
The finite differences \eqref{eq:data_finite_difference} are evaluated by using three different boundary value combinations
\[
 \eps(f_1+f_2), \ \  \eps f_1 \ \text{ and } \  \eps f_2,
\]
where %$f_1=(f_1^L,f_1^R)$ with~\f{We only have one $f=(f^L,f^R)$ here. Let's be more specific.} 
$f_1$ and $f_2$ are as in \eqref{eq:choice_of_fs}.
%\begin{equation}\label{eq:boundaryvalues}
%\begin{split}
%f_1(a,t_j) &= H^\tau\big((a-x_0)-(t_j-t_0)\big), \\ %\quad\text{and}\\
%f_2(b,t_j) &= H^\tau\big((b-x_0)+(t_j-t_0)\big).
%\end{split}
%\end{equation}
%
The measurement data (here for simplicity written without noise) is then obtained as the finite differences 
\begin{equation}\label{eq:bvs}
\begin{split}
F^L(t_j) &= \frac{1}{\eps^2} \Big( \p_\nu u_{\eps(f_1+f_2)}(a,t_j) - \p_\nu u_{\eps f_1}(a,t_j) - \p_\nu u_{\eps f_2}(a,t_j)   \Big),\\
F^R(t_j) &= \frac{1}{\eps^2} \Big( \p_\nu u_{\eps(f_1+f_2)}(b,t_j) - \p_\nu u_{\eps f_1}(b,t_j) - \p_\nu u_{\eps f_2}(b,t_j)   \Big),
\end{split}
\end{equation}
where $\eps>0$ is chosen experimentally. Here $\p_\nu$ stands for the normal derivatives as in \eqref{eq:Neumann_derivatives}.

\smallskip
\noindent\textbf{Approach (2): Obtaining $\p_{\eps_1\eps_2}^2|_{\eps_1=\eps_2=0}\widetilde\Lambda$ via regularization}
Let $d>0$ be a small number.
We divide the square $[-d,d]^2$ into $N_\eps\in \N$ equally spaced grid points $(\epsilon_i,\epsilon_j)$, and such that the point $(0,0)$ is contained in the grid.
The noisy measurements of the DN map are then evaluated at the grid points as
$$
\widetilde\Lambda_{i,j}:=\Lambda_q(\epsilon_if_1+\epsilon_jf_2) + \E_{i,j},
$$
where $f_1$ and $f_2$ are as in \eqref{eq:choice_of_fs}.
Let $A$ denote the $N_\eps\times N_\eps$ anti-differentiation matrix \eqref{eq:discrete_derivative}. 
The regularized mixed derivative $F^L$ corresponding to measurements on the left side, $x=a$, of the lateral boundary is obtained by solving
\begin{equation}\label{eq:min_tikhonov}
AF^L = \widetilde\Lambda|_{x=a}.
\end{equation}
This is then an approximation for $\p_{\eps_1\eps_2}^2\widetilde\Lambda(\eps_1f_1+\eps_2f_2)$ at  $x=a$. As explained earlier we consider \eqref{eq:min_tikhonov} the minimization problem \eqref{eq:tikhonov} and use Tikhonov regularization. The regularized mixed derivative $F^R$ corresponding to measurements on $x=b$ is obtained similarly. % each boundary point $(a,t_j)$ and $(b,t_j)$, $j=1,\ldots,N_t+1$ 
Then 
\[
 (F^L,F^R)|_{\eps_1=0, \s \eps_2=0} \text{ is approximation of } \p_{\eps_1\eps_2}^2|_{\eps_1=\eps_2=0}\Lambda_q.
\]

In either of the approaches \textbf{(1)} or \textbf{(2)}, let us also define
\begin{align*}
f_{0}^L(t_j) &= H^{\tau=1}((a-x_0)-(t_j-t_0))\\
f_{0}^R(t_j) &= H^{\tau=1}((b-x_0)+(t_j-t_0)).
\end{align*}
Then $f_{0}^L$ and $f_{0}^R$ correspond to the discretization of the boundary value $\psi$ of the measurement function $v_0$ on the left and right side of the lateral boundary respectively.
The integral in formula \eqref{est:2D_noise},
\[
 -\frac{1}{2\pi}\int_{\Sigma} \psi\,\p_{\eps_1,\eps_2}^2\Big|_{\eps_1=\eps_2=0}\widetilde\Lambda(\eps_1f_1 +\eps_2f_2) \d\Sigma 
\]
is then our numerical solution to the inverse problem at $(x_0,t_0)$. The mixed derivative $\p_{\eps_1,\eps_2}^2|_{\eps_1=\eps_2=0}\widetilde\Lambda(\eps_1f_1 +\eps_2f_2)$ is numerically computed either according to approach \textbf{(1)} or \textbf{(2)}. The integral is evaluated numerically as a Riemann sum and we set

\smallskip
\noindent\framebox{
    \parbox{\textwidth}{%
    \smallskip
    \begin{equation}\label{eq:numerical_integral_id}
q^\mathrm{numerical}(x_0,t_0) := -\frac{1}{2\pi}\sum_{j=1}^{N_t+1} \Big(f_{0}^L(t_j)F^L(t_j) + f_{0}^R(t_j)F^R(t_j)\Big) \Delta t.
\end{equation}
%     The points $x_1\in M_1$ and $x_2\in M_2$ are to be identified if and only if for every $f\in C^\infty(\p M)$ the harmonic extensions $u_f^1$ and $u_f^2$ of $f$ to $(M_1, g_1)$ and $(M_2, g_2)$ satisfy
%     \begin{equation}\label{boxref}
%     u_f^1(x_1)=u_f^2(x_2).
%     \end{equation}
    }%
}
\smallskip

% \begin{equation}\label{eq:numerical_integral_id}
% q^\mathrm{numerical}(x_0,t_0) := \frac{1}{2\pi}\sum_{j=1}^{N_t+1} \Big(f_{0,a}(t_j)F_{a}(t_j) + f_{0,b}(t_j)F_{b}(t_j)\Big) \Delta t.
% \end{equation}
%
%This integral $q^\mathrm{numerical}(x_0,t_0)$ will be our numerical solution to the inverse problem.

\begin{remark}
When studying inverse problems with synthetic data one should keep in mind the so-called inverse crime:
If the same model or grid is used for both the forward problem and the inversion, the results of the inversion can in some cases be better than they would be in reality.
In our considerations we do not commit an inverse crime, since the reconstruction grid is independent from the grid used in the forward model.
\end{remark}
%\begin{remark} 
%While the effect of noisy measurements in the naive finite difference approach to the inverse problem is already considered in \cite{LLPT20}, we can actually say more.
%%
%Consider noisy (continuous) measurement of the DN map
%$$
%\tilde \Lambda_q = \Lambda_q + \E(t),
%$$
%where $\E$ is a random variable.
%%
%Then the integral \eqref{eq:integral_identity_derivitve} we are interested in becomes
%\begin{multline*}
%\int_0^T f_0(t) D_{\eps_1\eps_2}^2\Big|_{\eps_1=\eps_2=0}\tilde\Lambda_q(\eps f_1+ \eps f_2) \d t
%= 
%\int_0^T f_0(t) D_{\eps_1\eps_2}^2\Big|_{\eps_1=\eps_2=0}\Lambda_q(\eps f_1+ \eps f_2) \d t\\ + D_{\eps_1\eps_2}^2\Big|_{\eps_1=\eps_2=0} \int_0^T f_0(t) \E(t) \d t.
%\end{multline*}
%Considering $f_0$ as a random variable independent of $\E$, the latter term becomes
%$$
%\int_0^T f_0(t) \E(t) \d t = \mathbb{E}(f_0(t)) \mathbb{E}(\E(t)).
%$$
%In particular, when the noise $\E$ has zero expected value, the contribution of the noise vanishes from the integrals.
%%
%Heuristically, when the measurement $\tilde \Lambda_q$ contains additive noise with zero expected value, the integral identity averages out the noise.
%\end{remark}

\subsection{Choosing  the parameters $\tau$ and $\eps$}\label{sec:tau_and_eps}

For the inverse problem we need to choose the parameters $\tau>0$ and $\eps>0$ for the boundary values. We remark that the work \cite{LLPT20}, given a bound for the noise $\mathcal{E}: H^{s+1}(\Sigma)\to H^r(\Sigma)$, with $s\in \N$ and $r\leq s$, computed theoretical values for $\tau$ and $\eps$ that can be used to achieve H\"older stability of the recovery in the inverse problem. In our numerical simulation we however encounter problems using these theoretically values. (Consequently, we expect that the  theoretical values for $\tau$ and $\eps$ used in \cite{LLPT20} are not optimal in general.)

\begin{figure}[ht!]
\centering
\includegraphics[height=4.3cm]{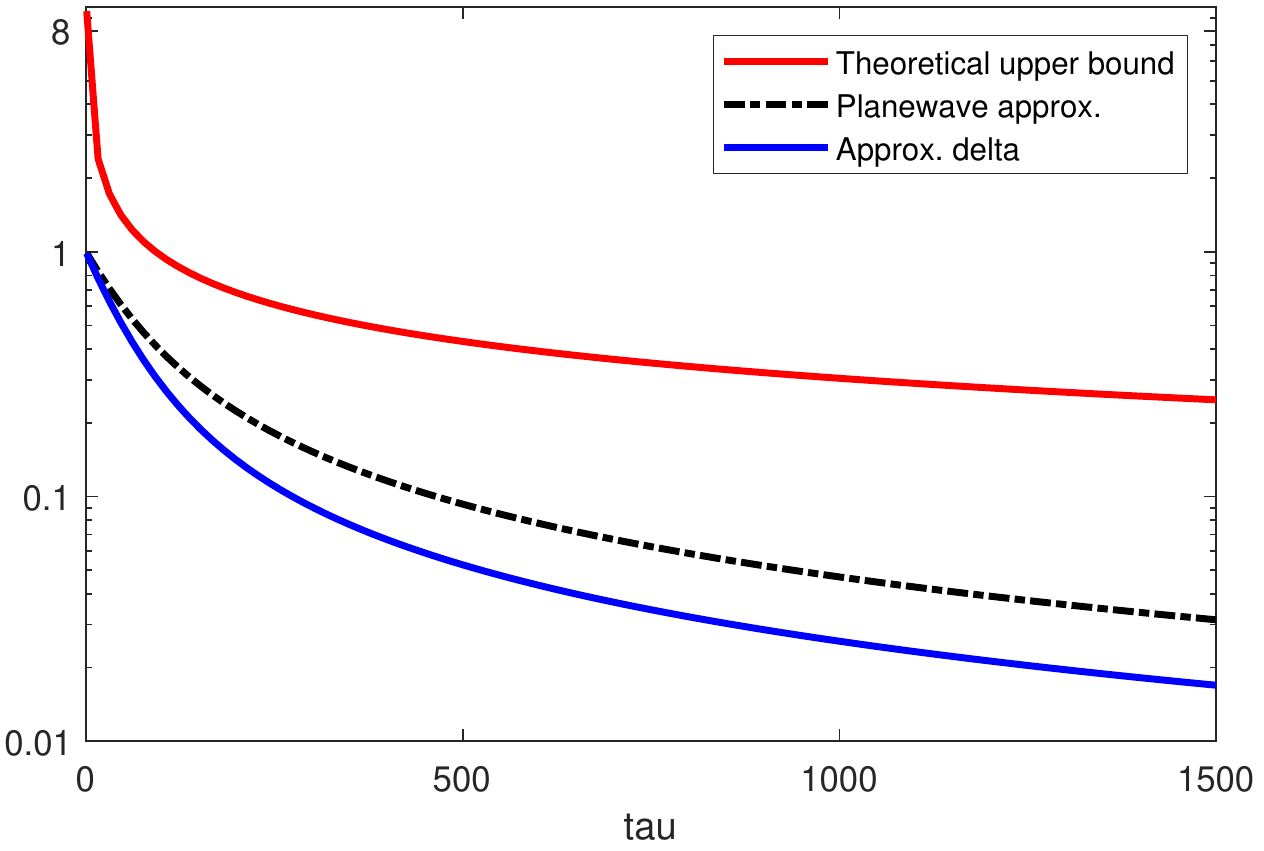}\hfill
\includegraphics[height=4.3cm]{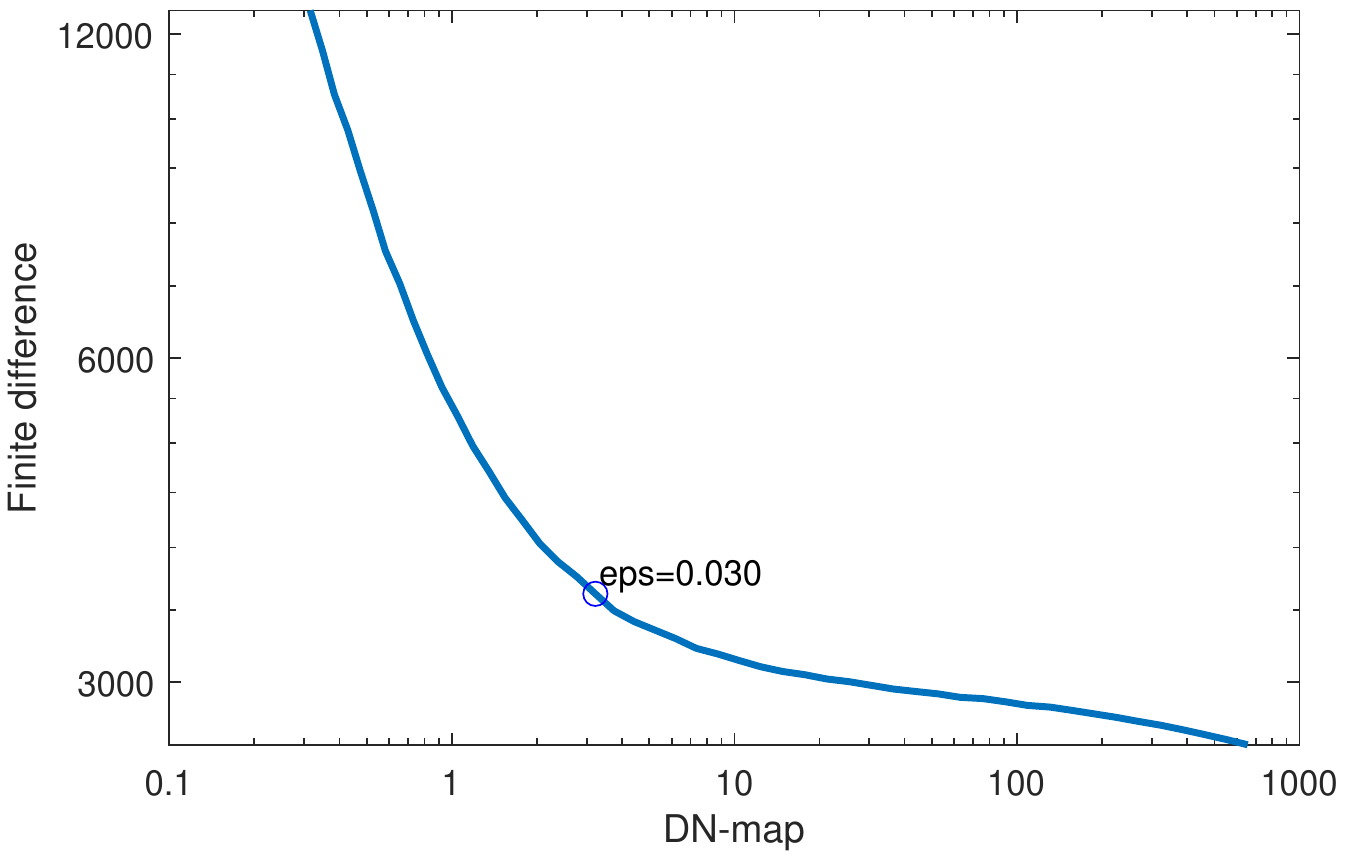}
\caption{Left: Convergence rate for the integral in Lemma~\ref{lemma:tau} with respect to $\tau\in [0,1500]$.
We compare the absolute errors between the integral approximations in sup-norm. Here the red curve depicts the theoretical upper bound $O(\tau^{-1/2})$ given by Lemma~\ref{lemma:tau}, the blue curve corresponds to the precise integral of Lemma~\ref{lemma:tau}  and the black dot-dash curve depicts the integral of \eqref{eq:tau2}.
The difference in convergence rates in both pictures is explained by the test function $\chi$ included in the functions $H_1^\tau$ and $H_2^\tau$ in the integral \eqref{eq:tau2}.
Right: The L-curve obtained by comparing the DN maps against the finite differences \eqref{eq:bvs} in $\log\log$-scale. In this case, the ``corner" of the curve is located approximatively at $\eps=0.03$.
 \label{fig:tau_and_eps}}
\end{figure}

We choose the parameters $\tau$ and $\eps$ as follows. The parameter $\tau>0$ is chosen by using the following lemma, whose proof can be found from \cite{LLPT20}.
\begin{lemma}\label{lemma:tau}
Assume $q$ is a Lipschitz function with Lipschitz constant $L>0$ and let $\tau>0$. The following estimate 
\[
\left| q(x_0,t_0)-  \frac{\tau}{\pi} \int_{\mathbb{R}^2}q(x,t)\e^{-\tau((x-x_0)^2 + (t-t_0)^2)} \d x\, \d t \right| \leq \frac{\sqrt{\pi}}{2} L\s \tau^{-1/2}
\]
holds true for all $(x_0, t_0)\in \mathbb{R}^2$. In particular, the integral on the left converges uniformly to $q$ when $\tau \to \infty$.
\end{lemma}
We evaluated the integral in Lemma~\ref{lemma:tau} numerically via Simpson's quadrature rule for the function 
$$q(x,t) = \varphi(5\sqrt{x^2+(t-1.5)^2})$$
 with various different values for $\tau$.
The absolute errors are depicted in Figure~\ref{fig:tau_and_eps}. 
Evidently larger values of $\tau>0$ result in better approximations.
%
%However choosing too large $\tau$ presents some problems in the numerical simulation of the wave equation \eqref{eq:intro_wave-eq}.
%
In our simulations, we used the value $\tau=700$. % \tbl{which seems to be enough.}
%
%\begin{figure}[ht!]
%\centering
%\includegraphics[width = 0.45\textwidth]{img/approx_pointwise_tau_EX1}\hfill
%\includegraphics[width = 0.45\textwidth]{img/approx_supnorm_tau_EX1}
%\caption{Convergence rates for the integral in Lemma~\ref{lemma:tau} with respect to $\tau\in [0,1500]$. Left: Comparison of the convergence of the integral in Lemma~\ref{lemma:tau} and the integral \eqref{eq:tau2} to the true value $q(x_0,t_0)=1$. 
%Right: Comparison of the absolute errors between the integral approximations in sup-norm. Here the red curve depicts the theoretical upper bound $O(\tau^{-1/2})$ given by Lemma~\ref{lemma:tau}, the blue curve corresponds to the precise integral of Lemma~\ref{lemma:tau} and the black dot-dash curve depicts the integral of \eqref{eq:tau2}.
%The difference in convergence rates in both pictures is explained by the test function $\chi$ included in the functions $H_1^\tau$ and $H_2^\tau$ in the integral \eqref{eq:tau2}.
% \label{fig:tau}}
%\end{figure}

To numerically compute the derivative $\p_{\eps_1\eps_2}^2 \Lambda_q(\eps_1f_1+\eps_2f_2)$ using finite differences, one has to also make a choice of $\eps$. This is in accordance to the formula \eqref{eq:bvs}.
Our choice is based on the following heuristic method, often called the L-curve method in the literature \cite{Hansen1,Hansen2,MS2012}.
We evaluated the synthetic DN map with multiple values for $\eps\in[0.001,1]$ with logarithmically spaced points.
Using these values for $\eps>0$, we then evaluated the corresponding DN maps $\Lambda_q(\eps f_1+\eps f_2)$, $\Lambda_q(\eps f_1)$ and $\Lambda_q(\eps f_2)$.
The noisy measurements are simulated by adding Gaussian noise with standard deviation equal to 
%\f{Why is size of the noise is chosen this way? ---Added the signal-to-noise ratios here. -TT}
\begin{equation}\label{eq:noise}
\sigma\cdot \mathrm{mean}_j(|\p_\nu u(a,t_j)|+|\p_\nu u(b,t_j)|)
\end{equation}
at the boundary to the synthetic DN map.
For this we used Matlab's \textsc{randn} function and the values $\sigma = 0.005,\s 0.01,\s0.015,\s0.02,\s0.025,\s0.03$ (signal-to-noise ratio $42.2$--$57.8\,\mathrm{dB}$).
Now, plotting the $L^2$-norms of the boundary data: 
$$
\Vert \Lambda_q(\eps f_1+\eps f_2) \Vert_{L^2(\Sigma)} +
\Vert \Lambda_q(\eps f_1) \Vert_{L^2(\Sigma)} +
\Vert \Lambda_q(\eps f_2) \Vert_{L^2(\Sigma)}
$$
against the $L^2$-norms of the finite differences $F^L(t_j)$ and $F^R(t_j)$ defined in \eqref{eq:bvs} in $\log\log$-scale, we get the L-shaped curve shown in Figure~\ref{fig:tau_and_eps}.
Heuristically, the reason for the shape of this curve is that when $\eps\to 0$ two things happen: Firstly, the DN map of $\eps f_1+\eps f_2$ etc. tend to zero.
Secondly, due to the presence of noise in the data, the finite differences $F^L$ and $F^R$ become numerically unstable, increasing their norm.
Therefore, one can try to balance between the size of the DN map of $\eps f_1+\eps f_2$ etc. and the error one makes using the finite differences.
Our choice of $\eps$ %\f{Teemu, do we use this $\eps$ in the examples? ---No, the $\eps$ we choose depends on the example. The $\eps$ is calculated by using the L-curve method to each problem separately. I removed the specific value of $\eps$. -TT} 
lies at the ``corner" of the L-curve, where both the norm of the measured (noisy) DN map and the norm of the finite differences are small. 

\begin{remark}
Note that for approach \textbf{(2)} one does not make a choice of a specific $\eps$. Instead, one uses several values of $\eps_1,\eps_2$ on some interval containing $\eps_1=\eps_2=0$ and then finds the derivative $\p_{\eps_1\eps_2}^2 \widetilde\Lambda(\eps_1f_1+\eps_2f_2)$ as a solution to a minimization problem.
\end{remark}
%
%\begin{figure}[ht!]
%\centering
%\includegraphics[width = 0.8\textwidth]{img/Lcurve_EX1sigma5-0e-03}
%\caption{The L-curve obtained by comparing the DN maps against the finite differences \eqref{eq:bvs} in $\log\log$-scale. In this case, the ``corner" of the curve is located approximately at $\eps=0.03$. \label{fig:Lcurve}}
%\end{figure}

%\newpage
\section{Numerical examples}\label{sec:examples}
We used the parameters $N_x = 1000$ and $N_t = 20000$ to solve the forward problem in the domain $[-0.5,0.5]\times[0,3]$.
These choices of $N_x$ and $N_t$ yield the CFL number $c=0.15$. 
For the inverse problem, we considered a $N\times N = 20\times 20$ equally spaced reconstruction grid of the domain $[-0.4,0.4]\times[1.3,1.7]\subset [-0.5,0.5]\times[0,3]$. As the parameter $\tau$ we use $\tau=700$.
\\
\noindent
\textbf{Approach 1: Finite differences} When using finite differences to approximate $\p_{\eps_1\eps_1}^2|_{\eps_1=\eps_2=0}\widetilde\Lambda(\eps_1f_1+\eps_2f_2)$, we used $\eps \in [0.01, 0.2]$ chosen by the L-curve method explained in Section~\ref{sec:tau_and_eps}.
\\
\noindent
\textbf{Approach 2: Regularization} When using the regularization approach to approximate $\p_{\eps_1\eps_1}^2|_{\eps_1=\eps_2=0}\widetilde\Lambda(\eps_1f_1+\eps_2f_2)$ discussed in Sections~\ref{sec:Dregularization} and~\ref{sec:IP}, we used $15^2=225$ equally spaced values $(\epsilon_m,\epsilon_n) \in [-0.15,0.15]^2$.

We remark that the reconstructions obtained by these two approaches are visually indistinguishable.
Approach \textbf{(1)} has the benefit of being simple to implement, but could become unstable depending on the noise. Approach \textbf{(2)} is more involved, but produces more precise mixed derivatives of the noisy measurements.
For the noisy measurements we used Gaussian noise with standard deviations given in~\eqref{eq:noise}. These noise levels correspond to signal-to-noise ratios in the range of $42.2$--$57.8\,\mathrm{dB}$ (calculated using Matlab's \textsc{snr} function from the Signal Processing Package).
We mention without details that we got similar results by using uniformly distributed noise with zero mean instead of Gaussian noise.

Let now $\chi_A$ be the characteristic function of a set $A$, given by
$$
\chi_A(x)=
\begin{cases}
1,& x\in A,\\
0,& x\not\in A.
\end{cases}
$$
We test the reconstruction of a potential function from the measurements of noisy DN maps on the following examples:
\begin{enumerate}
\item $q(x,t) = \varphi(5\sqrt{x^2+(t-1.5)^2})$,
\item $q(x,t) = (\sin(4\pi x) + \sin(8\pi t)) \varphi(5\sqrt{x^2+(t-1.5)^2})$,
\item $q(x,t) = \chi_\mathrm{A=rectangle}(x,t)$, 
\item $q(x,t) = \varphi\Big(5\big(x-0.15\sin(\frac{\pi}{2}\cos(3\pi t))\big)\Big)$, (smooth bump oscillating in time).
\item $q(x,t) = \varphi(10x/3)\sin(x)$,
\end{enumerate}
Four of the five examples are intentionally chosen to be \emph{time-dependent}. This is to study the contrast to inverse problems for the linear Shr\"odinger wave equation. In these problems time-dependence of the potential can be problematic in both theoretical and numerical reconstructions. %\f{Matti, can we have a couple of references to demonstrate the problems for linear models?} 
The theoretical work \cite{LLPT20} does not consider cases where the potential $q$ is discontinuous. We however still test numerically such a case  in Ex. 3.

\begin{figure}[tp]
\centering
\includegraphics[width = 0.9\textwidth]{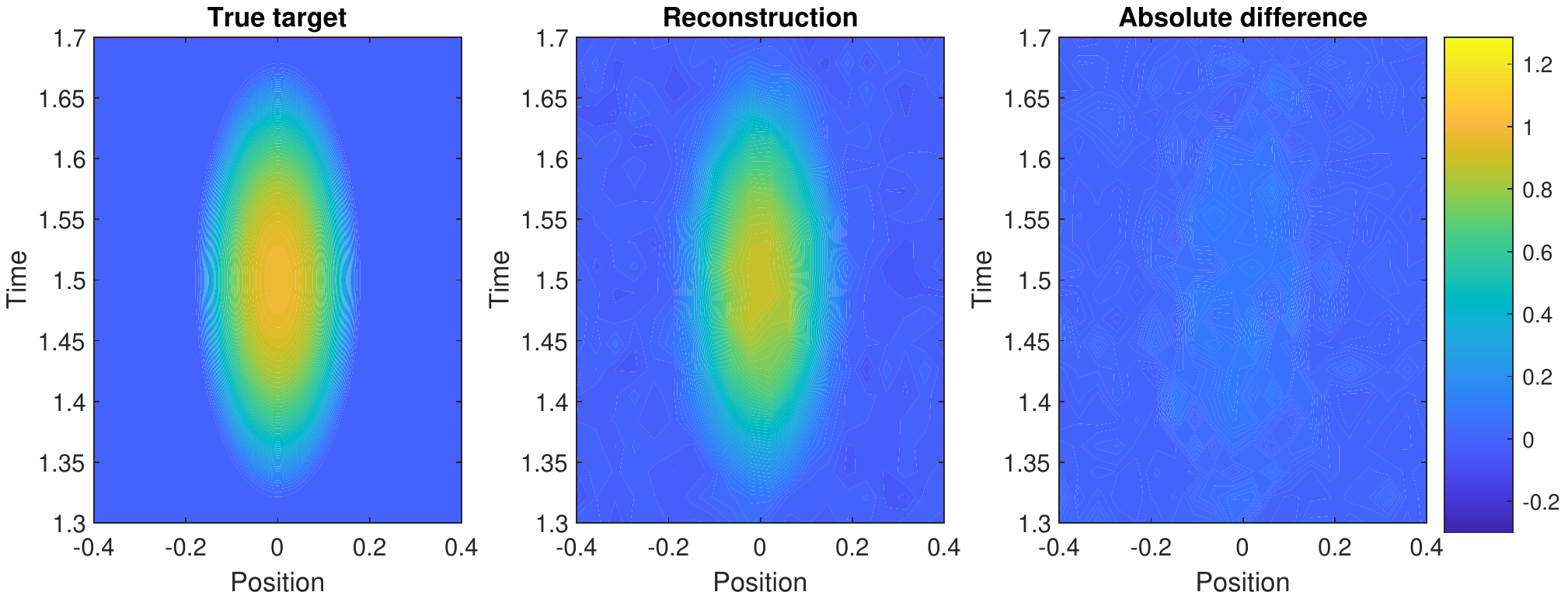}
\caption{Example 1: Comparison of the true target, supported on the time-interval, and its numerical reconstruction. \label{fig:EX1}}
\end{figure}

\begin{figure}[tp]
\centering
\includegraphics[width = 0.9\textwidth]{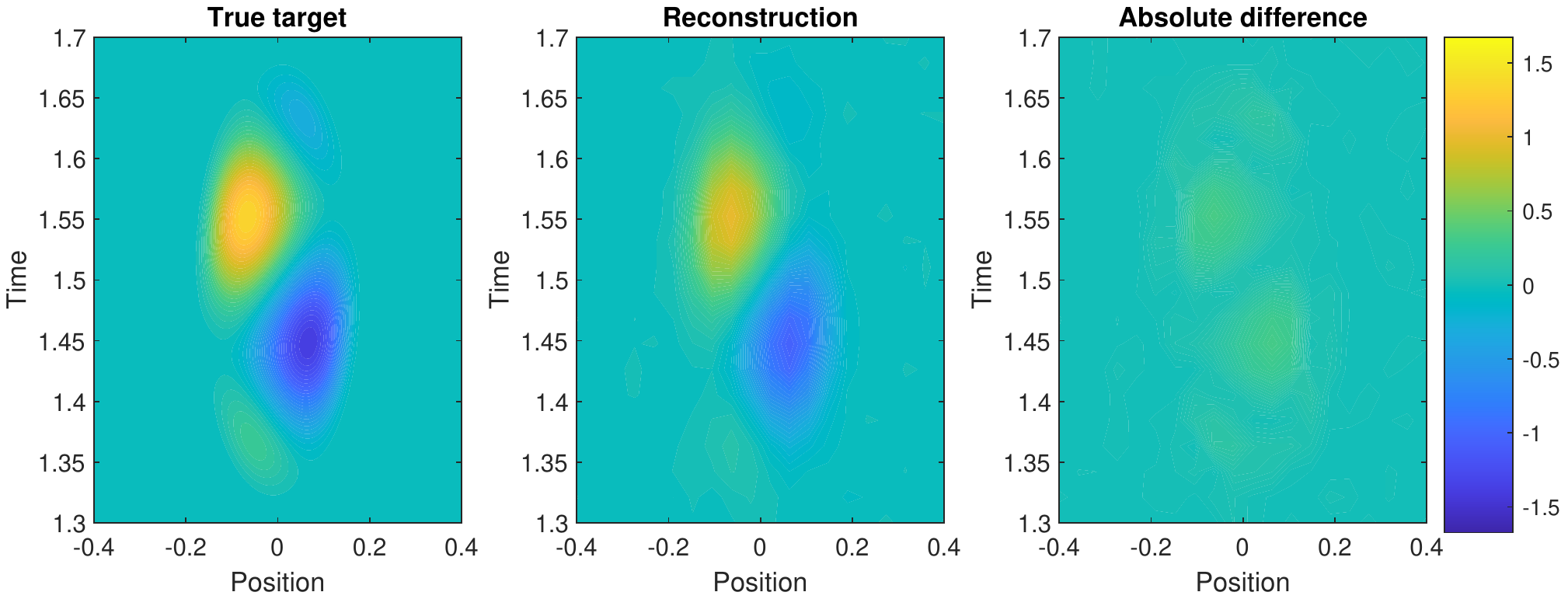}
\caption{Example 2: Comparison of the true target, which changes sign on the time-interval, and its numerical reconstruction. \label{fig:EX2}}
\end{figure}

\begin{figure}[tp]
\centering
\includegraphics[width = 0.9\textwidth]{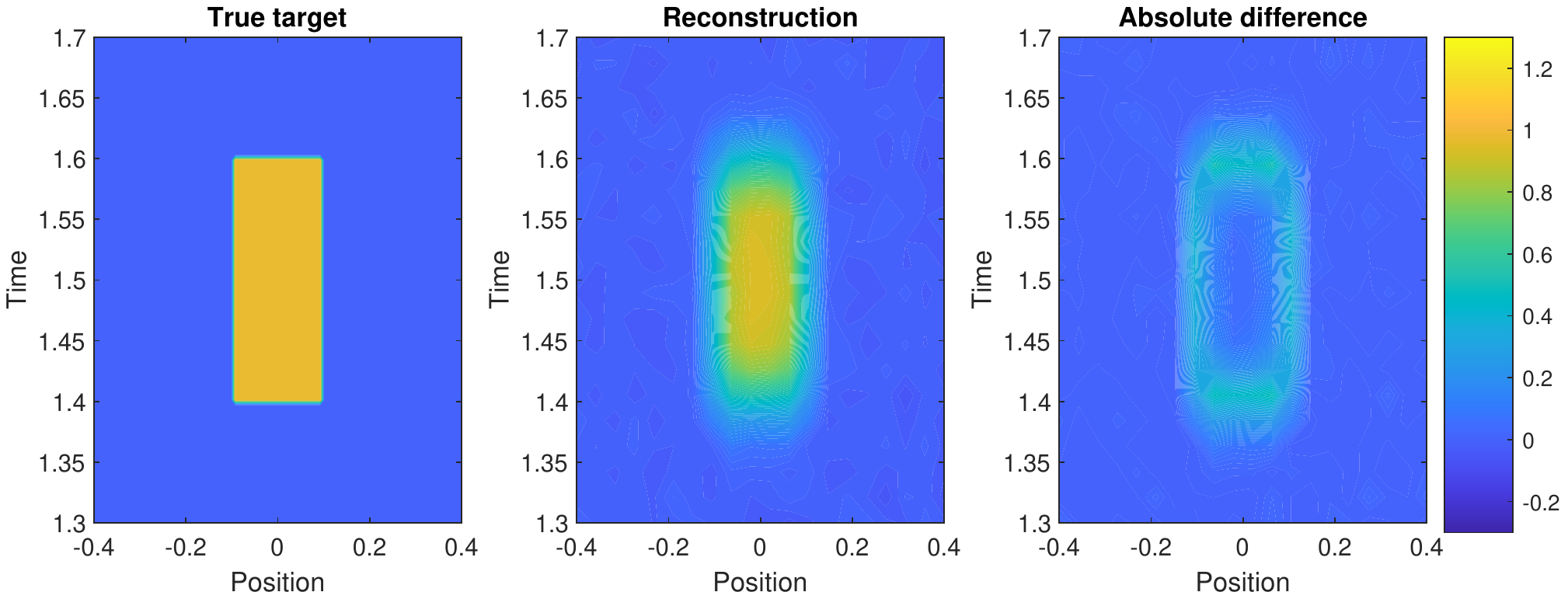}
\caption{Example 3: Comparison of the true target and its numerical reconstruction in the presence of discontinuities. The jump discontinuity of the target is blurred in the reconstruction. 
%This is caused \tbl{We know this for a fact, or speculate?} by the averaging property of the approximate identity integrals \eqref{eq:tau2}. 
\label{fig:EX3}}
\end{figure}

\begin{figure}[tp]
\centering
\includegraphics[width = 0.9\textwidth]{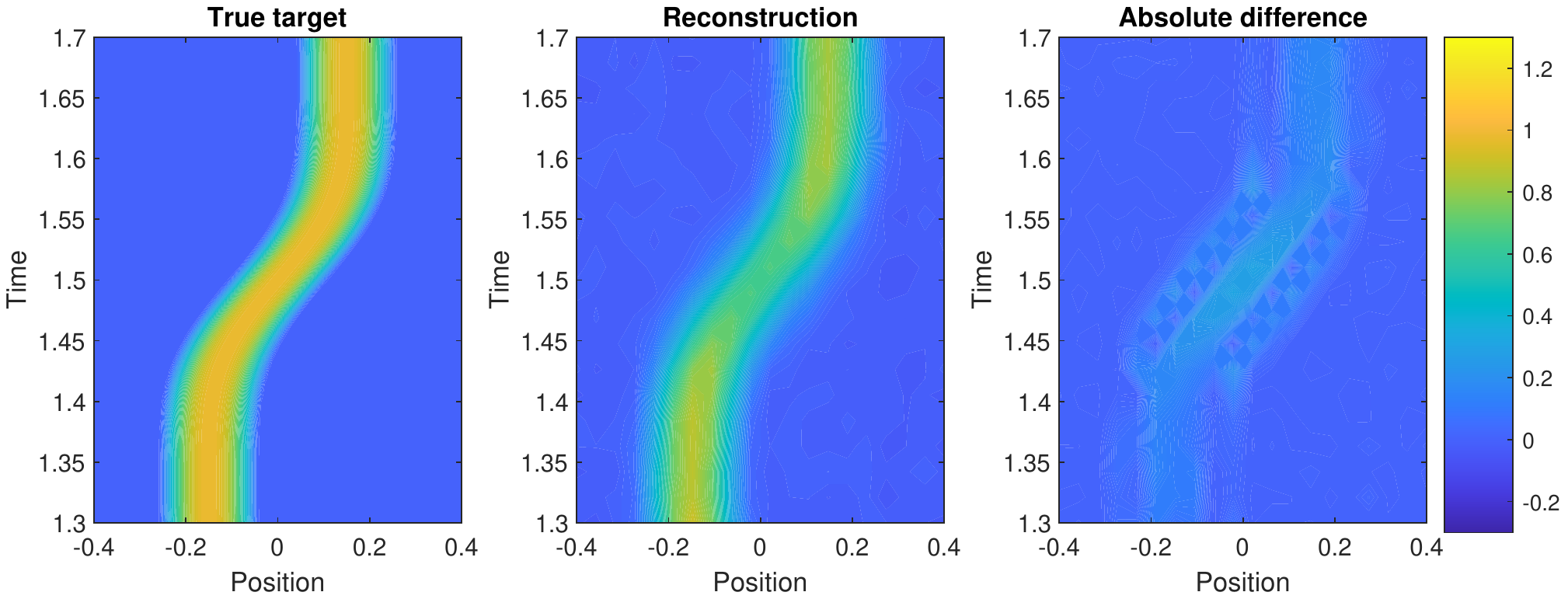}
\caption{Example 4: Comparison of the true non-compactly supported target target and its numerical reconstruction. The true target is a smooth bump-function, whose location moves in space along time.  \label{fig:EX4}}
\end{figure}

\begin{figure}[tp]
\centering
\includegraphics[width = 0.9\textwidth]{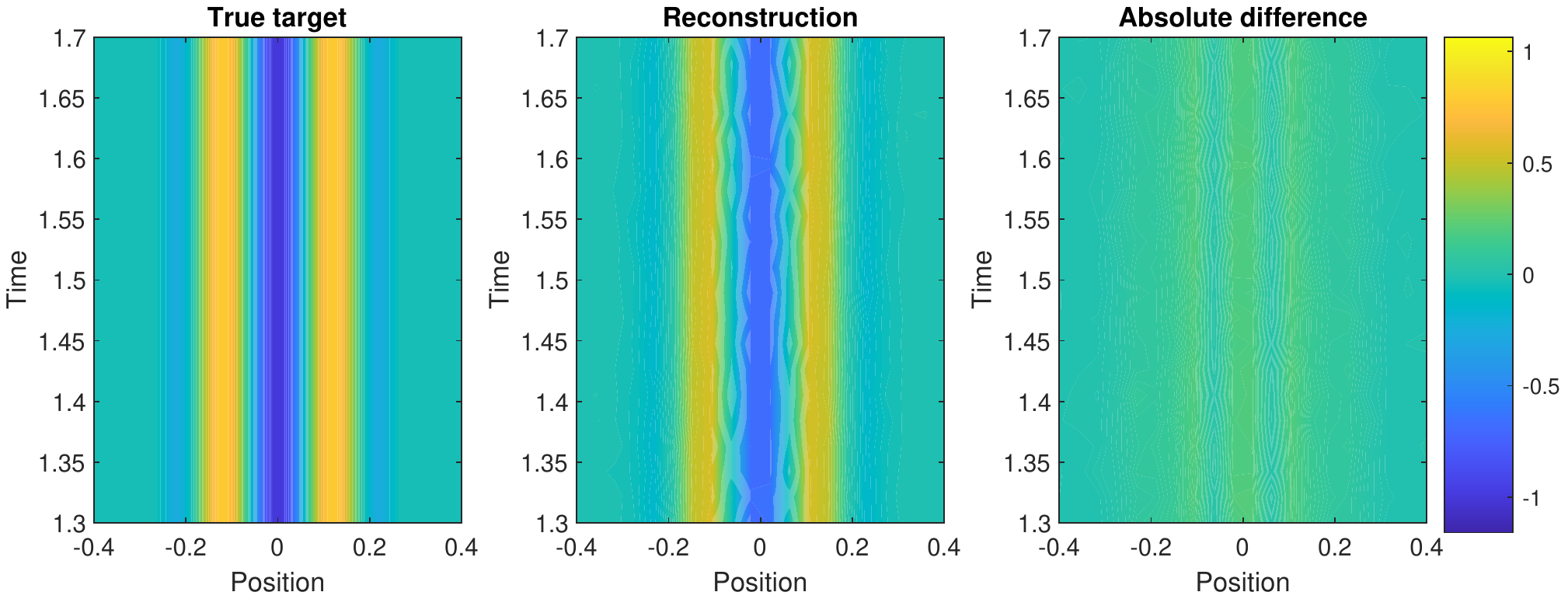}\\
\includegraphics[width = 0.9\textwidth]{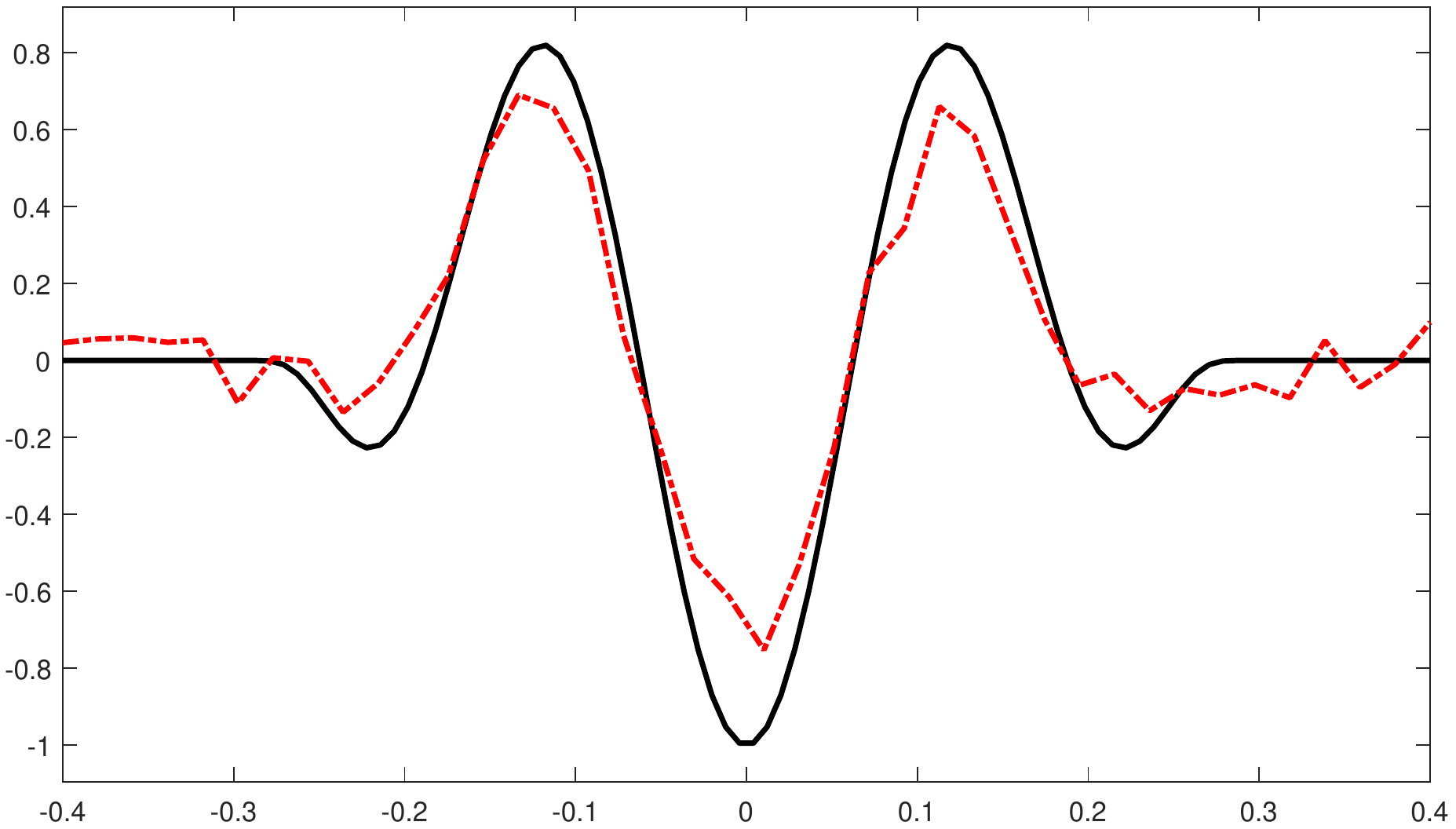}
\caption{Example 5: Comparison of the true  target and its numerical reconstruction in case of time-independent sign-changing target. Top: reconstruction of the potential function in time. Bottom: a cross-section of the true target and its reconstruction at $t_0=1.5$ in black and red dot-dash, respectively.  \label{fig:EX5}}
\end{figure}

Figures~\ref{fig:EX1}--\ref{fig:EX5} show the corresponding reconstructions for the five examples above obtained by evaluating the Riemann sums \eqref{eq:numerical_integral_id} over the reconstruction grid.
In all of the numerical reconstructions, the location and shape of the potential term $q$ is well visible to the eye.  In Table~\ref{tab:noise} we report the values for the absolute error 
$$
\eps_\mathrm{abs} = \max_{j}| q_j^\mathrm{numerical} - q_j^\mathrm{exact} |,
$$
the absolute relative error
$$
\eps_\mathrm{rel} = \frac{ \max_{j}| q_j^\mathrm{numerical} - q_j^\mathrm{exact} | }{ \max_{j}|q_j^\mathrm{exact} |}
$$
and the relative error in $L^2$-norm
$$
\eps_{L^2} = \frac{\sqrt{\sum_{j=1}^{N^2}(q_j^\mathrm{numerical} -q_j^\mathrm{exact} )^2 }}{\sqrt{\sum_{j=1}^{N^2}(q_j^\mathrm{exact} )^2 }}
$$
% and the root mean square error
% $$
% \eps_{\mathrm{RMSE}} = \sqrt{ \frac{1}{N^2}\sum_{j=1}^{N^2}(q_j^\mathrm{numerical} -q_j^\mathrm{exact} )^2 }
% $$
with respect to various noise levels $\sigma$.
%

%In Figure~\ref{fig:EX6} we also test reconstruction of
% $$
% q(x,t) = \chi_\mathrm{a}(\sqrt{x^2+t^2}) + \chi_\mathrm{b}(\sqrt{x^2+t^2}),
% $$
%where $\chi_\mathrm{a},\chi_\mathrm{b}$ are characteristic functions of two small intervals. This produces characteristic functions of two annuli. In this case the reconstruction has serious difficulties due to small contrast and discontinuities.

\begin{table}[tp]
\centering
\begin{tabular}{ |c|c|c|c|c|c|c|c|c|c| } 
\hline
& $\sigma$ & 0 & 0.005 & 0.01 & 0.015 & 0.02 & 0.025 & 0.03 \\
\hline 
\multirow{3}{2em}{Ex.1} 
	& $\eps_\mathrm{abs}$ & 0.088 & 0.115 & 0.130 & 0.135 & 0.144 & 0.160 & 0.166\\
	& $\eps_\mathrm{rel}$ & 0.089 & 0.116 & 0.132 & 0.137 & 0.146 & 0.163 & 0.169\\ 
	& $\eps_{L^2}$ 		  & 0.121 & 0.133 & 0.135 & 0.139 & 0.151 & 0.155 & 0.163\\
%& $\eps_{\mathrm{RMSE}} $ & 0.038 & 0.041 & 0.042 & 0.044 & 0.047 & 0.049 & 0.051\\
\hline
\multirow{3}{2em}{Ex.2} 
	& $\eps_\mathrm{abs}$ & 0.330 & 0.333 & 0.349 & 0.358 & 0.350 & 0.335 & 0.365\\
	& $\eps_\mathrm{rel}$ & 0.240 & 0.242 & 0.254 & 0.261 & 0.255 & 0.244 & 0.266\\ 
	& $\eps_{L^2}$ 		  & 0.250 & 0.253 & 0.259 & 0.261 & 0.258 & 0.255 & 0.260\\ 
%& $\eps_{\mathrm{RMSE}} $ & 0.074 & 0.075 & 0.077 & 0.078 & 0.077 & 0.076 & 0.077\\
\hline
\multirow{3}{2em}{Ex.3} 
	& $\eps_\mathrm{abs}$ & 0.494 & 0.519 & 0.505 & 0.508 & 0.520 & 0.508 & 0.505\\
	& $\eps_\mathrm{rel}$ & 0.494 & 0.519 & 0.505 & 0.508 & 0.520 & 0.508 & 0.505\\ 
	& $\eps_{L^2}$ 		  & 0.372 & 0.378 & 0.376 & 0.377 & 0.378 & 0.379 & 0.378\\ 
%& $\eps_{\mathrm{RMSE}} $ & 0.118 & 0.119 & 0.119 & 0.119 & 0.119 & 0.119 & 0.119\\
\hline
\multirow{3}{2em}{Ex.4} 
	& $\eps_\mathrm{abs}$ & 0.299 & 0.312 & 0.316 & 0.317 & 0.357 & 0.395 & 0.460\\
	& $\eps_\mathrm{rel}$ & 0.299 & 0.312 & 0.316 & 0.317 & 0.357 & 0.395 & 0.460\\ 
	& $\eps_{L^2}$ 		  & 0.206 & 0.216 & 0.238 & 0.256 & 0.275 & 0.292 & 0.315\\ 
%& $\eps_{\mathrm{RMSE}} $ & 0.077 & 0.081 & 0.089 & 0.096 & 0.103 & 0.109 & 0.118\\
\hline
\multirow{3}{2em}{Ex.5} 
	& $\eps_\mathrm{abs}$ & 0.198 & 0.216 & 0.221 & 0.222 & 0.228 & 0.223 & 0.231\\
	& $\eps_\mathrm{rel}$ & 0.260 & 0.283 & 0.290 & 0.290 & 0.298 & 0.292 & 0.302\\ 
	& $\eps_{L^2}$ 		  & 0.249 & 0.250 & 0.252 & 0.251 & 0.251 & 0.254 & 0.252\\ 
%& $\eps_{\mathrm{RMSE}} $ & 0.104 & 0.104 & 0.105 & 0.105 & 0.104 & 0.106 & 0.105\\
\hline
\end{tabular}
\caption{Comparison of various error quantities between the exact unknown $q^\mathrm{exact}$ and the numerical solution $q^\mathrm{numerical}$ under different noise levels $\sigma$. \label{tab:noise}}
\end{table}

%\begin{figure}[tp]
%\centering
%\includegraphics[width = 0.9\textwidth]{img/u_EX1}
%\caption{The unknown potential function (left) and the corresponding solution $u$ of \eqref{eq:intro_wave-eq} with parameters $(x_0,t_0)=(0,1.5)$ (right). This solution is the same as in Example 1 of Section~\ref{sec:examples} for the inverse problem. \label{fig:u}}
%\end{figure}


\begin{figure}[tp]
\centering
\includegraphics[width = 0.3\textwidth]{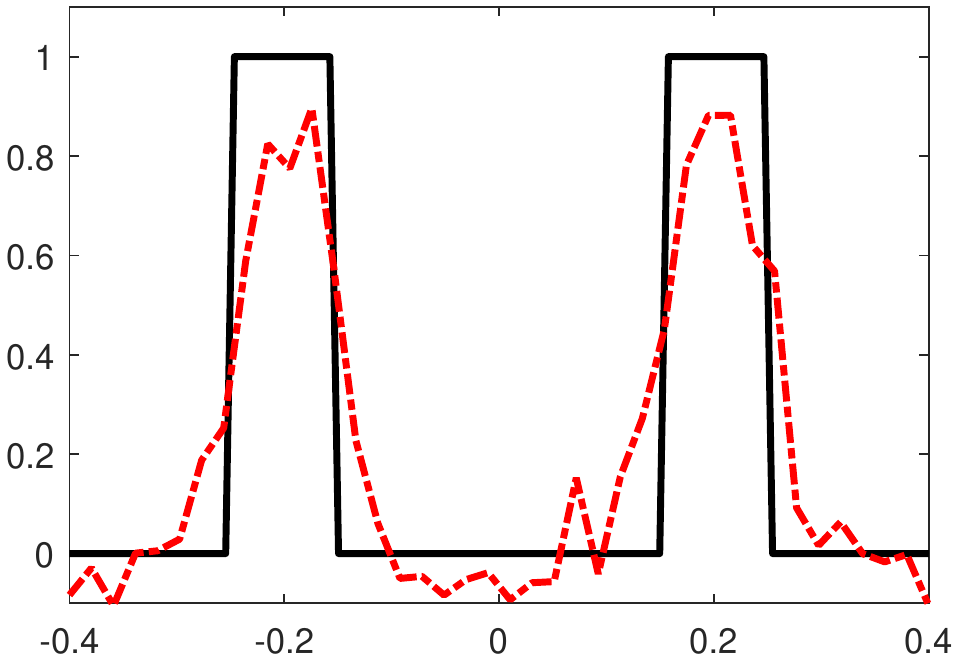}
\hfill
\includegraphics[width = 0.3\textwidth]{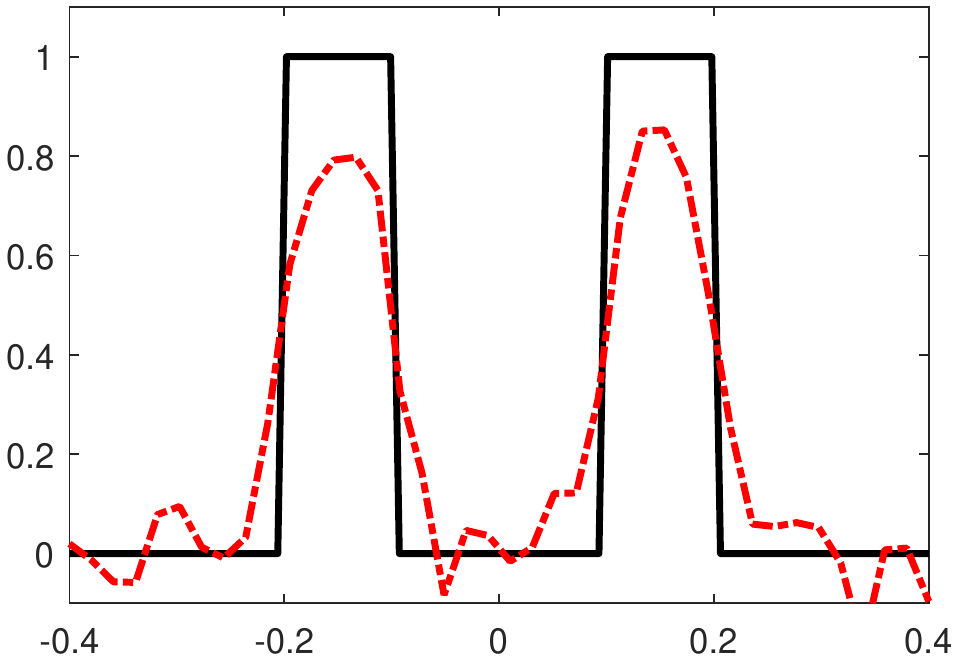}
\hfill
\includegraphics[width = 0.3\textwidth]{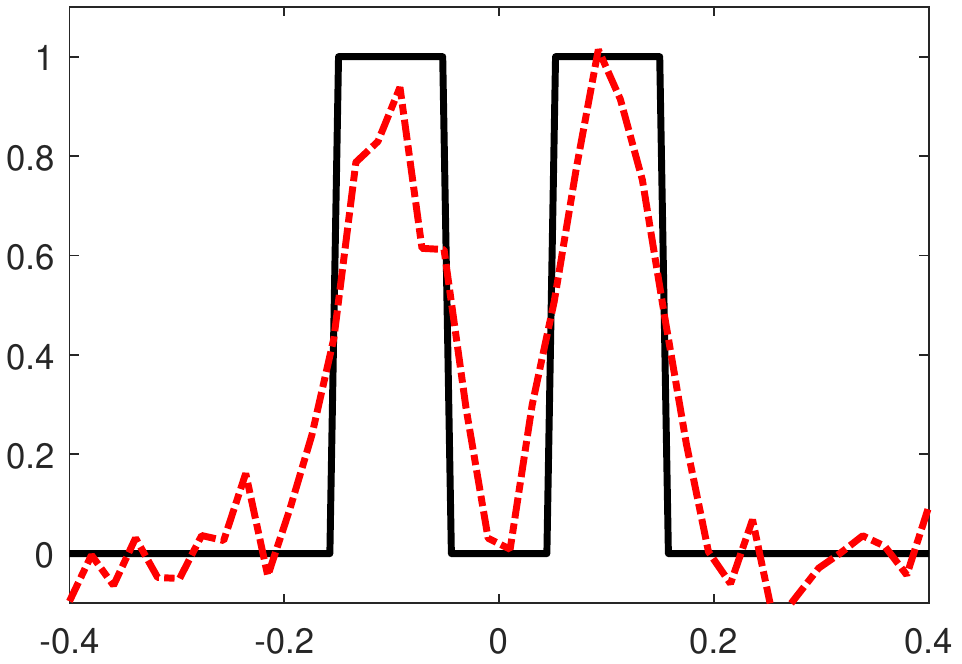}\\
\smallskip
\includegraphics[width = 0.3\textwidth]{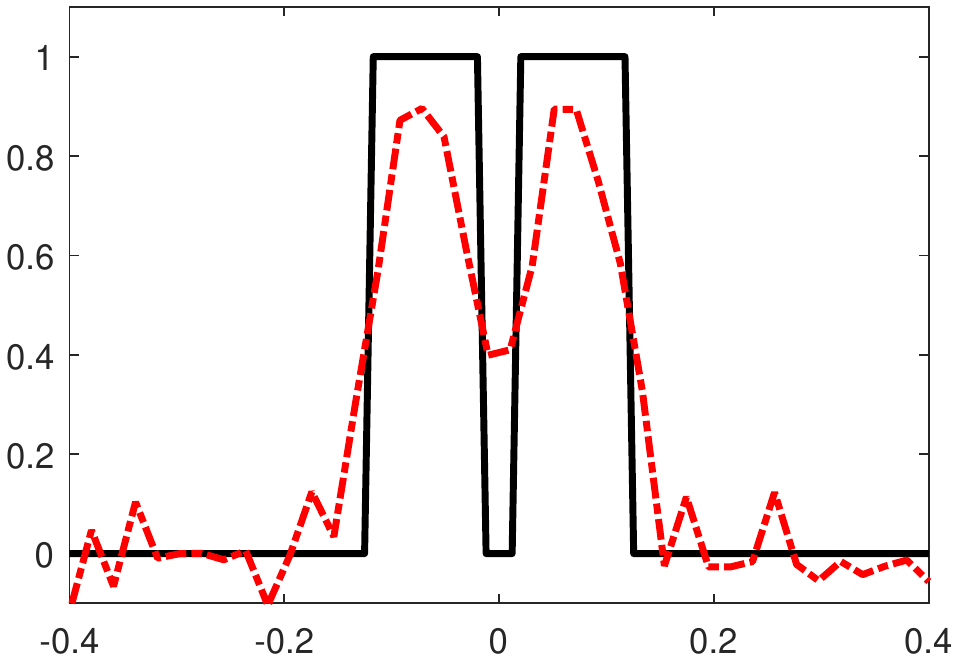}
\hfill
\includegraphics[width = 0.3\textwidth]{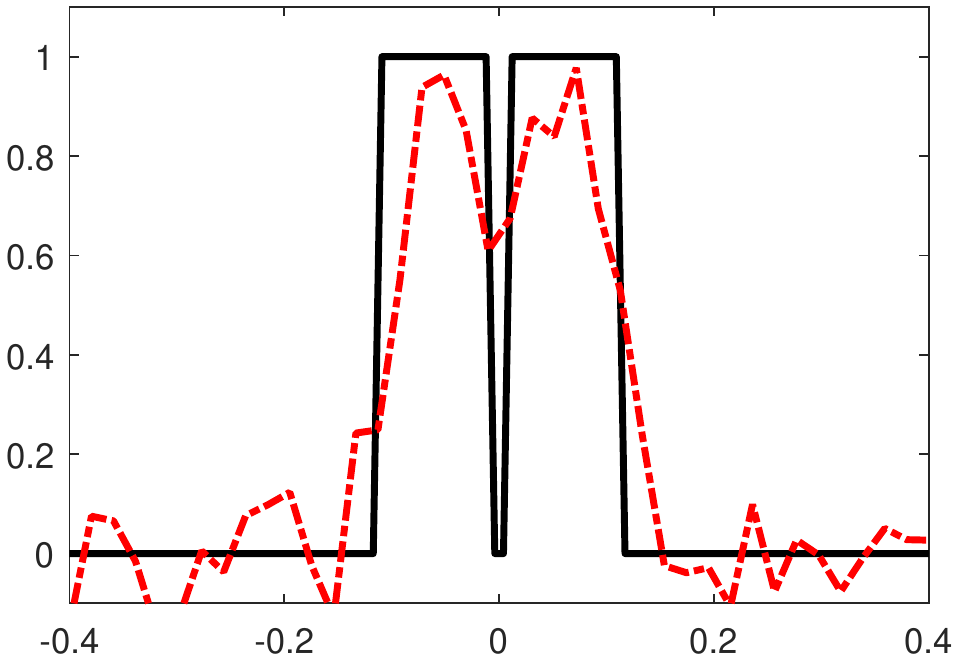}
\hfill
\includegraphics[width = 0.3\textwidth]{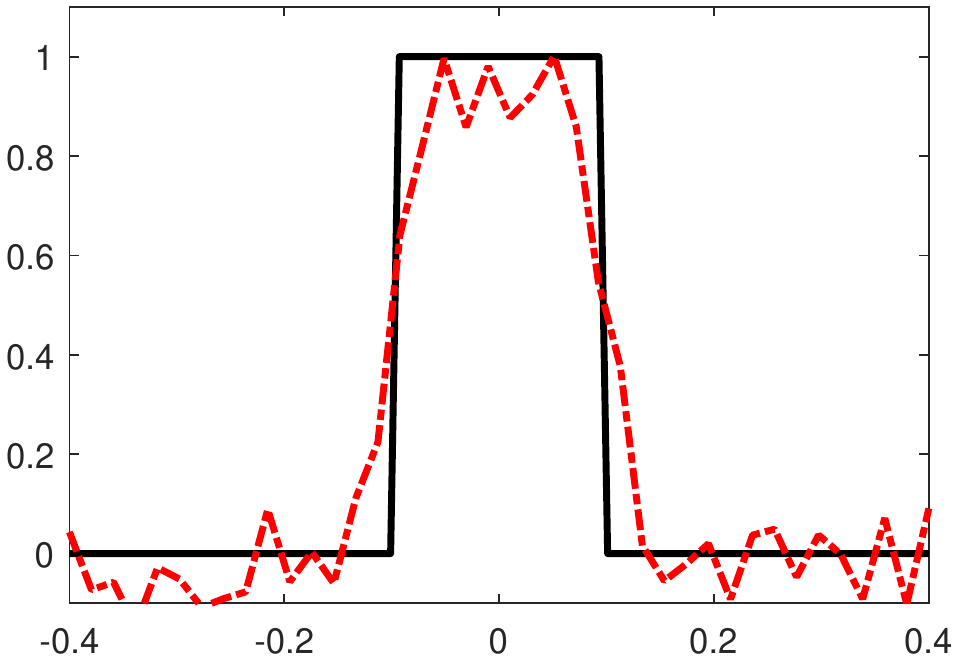}
\caption{Example 6: To demonstrate the resolution of our numerical method, we consider time-independent potentials $q(x)$ given by the characteristic functions of two small intervals of width $0.1$ at distances $0.4$, $0.3$, $0.2$, $0.14$, $0.12$ and $0.1$ as measured from the center points of the intervals. The cross-sections of the true targets and their reconstructions at fixed time $t_0=1.5$ are given in black and red dot-dash, respectively. Two sufficiently far apart target intervals are clearly distinguishable. However, we lose details when the targets are close to each other. There are two reasons for this. Firstly, our approximate plane waves have a finite width of approximately $0.1$ units, which causes loss of small details. Secondly, the potential functions considered are discontinuous and thus the jump discontinuity gets blurred (which is similar to Example 3 in Figure~\ref{fig:EX3}).
\label{fig:EX6}}
\end{figure}

%
%\begin{figure}
%\centering
%\includegraphics[width = 0.9\textwidth]{img/measurement_EX6_Sigma0-0e+00}\\
%\includegraphics[width = 0.9\textwidth]{img/measurement_EX6_Sigma5-0e-03}
%\caption{Example 6: Comparison of the true target and the numerical reconstruction in case of characteristic function of two annuli. Top figure shows the reconstruction when given unrealistic noiseless data. The bottom figure depicts the same reconstruction when $\sigma=0.005$. This is an example where the reconstruction is not so good. \label{fig:EX6}}
%\end{figure}

\section{Conclusions}
The methods discussed in this paper exploit the nonlinear nature of the wave equation \eqref{eq:intro_wave-eq}. This paper is a  numerical demonstration how a nonlinearity helps in inverse problems.

We studied an inverse problem for a one-dimensional nonlinear wave equation of the form $\square u + qu^2=0$ from a computational perspective. This was based on the theoretical reconstruction in \cite{LLPT20}.
The measurement data was the Dirichlet-to-Neumann map on the lateral boundary of the domain.
The synthetic DN map was evaluated by numerically solving \eqref{eq:intro_wave-eq} by using a finite difference scheme \eqref{eq:FDTD} and adding Gaussian noise.
An example of a solution to the forward problem and its numerical convergence were discussed.

For the inverse problem, we implemented two reconstruction algorithms corresponding to different ways of calculating mixed derivative of the DN map. The first one was based on finite difference approximation (as in \cite{LLPT20}) and the other was based on calculating the required mixed derivatives via a regularization method. The regularization method of calculating mixed derivatives numerically were discussed separately. 
The numerical approximation of the unknown potential was obtained as an integral against the numerical mixed derivative of the nonlinear DN map for specifically chosen boundary values.
We presented a heuristic method for choosing the required parameters $\tau$ and $\eps$ of the boundary values.
%

%The regularization method of calculating mixed derivatinves numerically were discussed separetely. 
%mixed derivatives approach to numerically calculating two dimensional second derivatives $\p_x\p_y$ was discussed.
%
Finally, multiple examples of reconstructions of potential functions were given, including smooth, discontinuous and time-dependent potentials.
The reconstruction was able to identify the location, shape, and size of the potential function.
%
%We report without details that if we drop the nonlinearity of \eqref{eq:intro_wave-eq} and consider instead the linear wave equation $\square u + au = 0$, then the reconstruction methods will produce a value zero for $a$. This is of course ex
%

%

\subsection*{Acknowledgements}

%M. L. was supported by Academy of Finland, grants 320113, 318990, and 312119.
L. P-M. and T. L. were supported by the Academy of Finland (Centre of Excellence in Inverse Modeling and Imaging, grant numbers 284715 and 309963) and by the European Research Council under Horizon 2020 (ERC CoG 770924). T. T  was partly supported by the Academy of Finland (Centre of Excellence in Inverse Modeling and Imaging, grant number 312119).

\noindent{\footnotesize E-mail addresses:\\
Matti Lassas: {matti.lassas@helsinki.fi}\\
Tony Liimatainen: {tony.liimatainen@helsinki.fi}\\
Leyter Potenciano-Machado: {leyter.m.potenciano@gmail.com}\\
Teemu Tyni: {teemu.tyni@utoronto.ca}
}


\begin{thebibliography}{99}


%
%
%\bibitem{AF03} Adams,~R. and Fournier,~J. Sobolev spaces, \emph{Academic Press}, 2003.
%
%\bibitem{AKKLT}
%Anderson,~M., Katsuda,~A., Kurylev,~Y., Lassas,~M., and Taylor,~M. \emph{Boundary regularity for the Ricci equation, geometric convergence, and Gel'fand's inverse boundary problem.} Invent. Math. \textbf{158},  261--321, 2004. 
%
%\bibitem{AZ17} Assylbekov,~Y.~M. and Zhou,~T. \emph{Direct and inverse problems for the nonlinear time-harmonic Maxwell equations in Kerr-type media}, to appear in J. Spectral Theory, arXiv:1709.07767, 2017
%
%\bibitem{AP} Astala,~K. and P\"aiv\"arinta,~L.  \emph{Calder{\'o}n's inverse conductivity problem in the plane}, Annals of Math., \textbf{163}, 265--299, 2006.
%
%
%\bibitem{Be87} Belishev,~M.I. \emph{An approach to multidimensional inverse problems for the wave equation}, Dokl. Akad. Nauk SSSR297.3, 524--527, 1987.
%
%\bibitem{BK92} Belishev,~M.I. and Kurylev,~Y. \emph{To the reconstruction of a Riemannian manifold via its spectral data (BC-method)}, Comm. PDE, \textbf{17}, 767--804, 1992.
%
%\bibitem{BL76} Bergh,~J. and L\"ofstrom,~J. Interpolation spaces: an introduction, \emph{Springer-Verlag}, Berlin, 1976
%
%\bibitem{BKL17} Bosi,~R., Kurylev,~Y., and Lassas,~M. \emph{Reconstruction and stability in Gel'fand's inverse interior spectral problem}, arXiv preprint arXiv:1702.07937, 2017.
\bibitem{BKLT20} Balehowsky,~T., Kujanp\"a\"a,~A., Lassas,~M., and Liimatainen,~T. \textit{An Inverse Problem for the Relativistic Boltzmann Equation}, (2020), arXiv:2011.09312.
%
%
%\bibitem{BK81}  Bukhgeim,~A. and Klibanov,~M. \emph{Uniqueness in the large of a class of multidimensional inverse problems}, Dokl. Akad. Nauk SSSR, \textbf{260}(2):269--272, 1981.
%
%\bibitem{Ca80} Calder{\'o}n,~P. \emph{On  an  inverse  boundary  value  problem}, Seminar on Numerical Analysis and its Applications to Continuum Physics, Soc. Brasileira de Matem{\'a}tica, Rio de Janeiro, 65--73, 1980.
%
%\bibitem{CNV19} C\^arstea,~C.~I., Nakamura,~G., and Vashisth,~M. \emph{Reconstruction for the coefficients of a quasilinear elliptic partial differential equation}, App. Math. Letters,  \textbf{98}, 121--127, 2019.
%
\bibitem{Ch2011} Chartrand,~R. \emph{Numerical Differentiation of Noisy, Nonsmooth Data}, ISRN Applied Mathematics, ID 164564, DOI:10.5402/2011/164564, 2011.

%\bibitem{Chen2019} Chen,~X., Lassas,~M., Oksanen,~L., and Paternain,~G. P. \emph{Detection of Hermitian connections in wave equations with cubic non-linearity}, arXiv:1902.05711, 2019.
%
%\bibitem{CLOP2020}
%Chen,~X., Lassas,~M., Oksanen,~L., and Paternain,~G. P. \emph{
%Inverse problem for the Yang-Mills equations.} arXiv:2005.12578, 2020.
%
%\bibitem{CB08} Choquet-Bruhat,~Y. General relativity and the Einstein equations. \emph{OUP Oxford}, 2008.
%
\bibitem{Cu1971} Cullum,~J. \emph{Numerical differentiation and regularization}, SIAM Journal on Numerical Analysis, vol. 8,
254–265, 1971.

%\bibitem{dH2018}
%de~Hoop,~M, Kepley,~P.,  and Oksanen,~L.,
%\emph{Recovery of a smooth metric via wave field and coordinate
%  transformation reconstruction},
%%\newblock \emph{Preprint arXiv:1710.02749}.
%SIAM J. Appl. Math., \textbf{78}, 1931--1953, 2018.
%
%
\bibitem{dH2019}de Hoop,~M., Uhlmann,~G. and Wang,~Y. \emph{Nonlinear responses from the interaction of two progressing waves at an interface}. Annales de l'Institut Henri Poincar\'e C, Analyse non lineaire, 36(2), 347--363, 2019.
%
\bibitem{dH2020} de Hoop,~M., Uhlmann,~G. and Wang,~Y. \emph{Nonlinear  interaction  of  waves  in  elastodynamics  and  an  inverse problem}. Mathematische Annalen, \textbf{376}(1-2), 765--795, 2020. 
%
%
%\bibitem{Esk} Eskin,~G. \emph{Inverse hyperbolic problems with time-dependent coefficients}, Commun. Partial Diff. Eqns., \textbf{32}(11), 1737--1758, 2007. 
%
%\bibitem{FIKO} Feizmohammadi,~A., Ilmavirta,~J., Kian,~Y. and Oksanen,~L. \emph{Recovery of time dependent coefficients from boundary data for hyperbolic equations}, To appear in J. of Spectral Theory, 2020.
%
%\bibitem{FeOk20} Feizmohammadi,~A. and Oksanen,~L. \emph{An  inverse  problem  for  a  semi-linear  elliptic  equation  in  Riemannian geometries}, J. Diff. Equations, \textbf{269}(6), 4683--4719, 2020.
%
%
%\bibitem{FO}
%Feizmohammadi,~A. and Oksanen, L. \emph{Recovery of zeroth order coefficients in non-linear wave equations}, To appear in J. Inst. Math. Jussieu, 2020.

\bibitem{FLO21} Feizmohammadi, A. and Lassas, M. and Oksanen, L., \textit{Inverse problems for nonlinear hyperbolic equations with disjoint sources and receivers}, Forum Math. Pi, \textbf{9}, (2021), \url{https://doi.org/10.1017/fmp.2021.11}. 

\bibitem{FO20} Feizmohammadi,~A. and Oksanen,~L. \emph{Recovery of zeroth order coefficients in non-linear wave equations}, To appear in J. Inst. Math. Jussieu, 2020.
%arXiv preprint arXiv:1903.12636 (2019).
%
\bibitem{FO19} Feizmohammadi,~A. and Oksanen,~L. \emph{An  inverse  problem  for  a  semi-linear  elliptic  equation  in  Riemannian geometries}, J. Diff. Equations, \textbf{269}(6), 4683--4719, 2020.


%
%\bibitem{GT13} Guillarmou,~C. and Tzou,~L. \emph{The Calder\'on inverse problem in two dimensions}, Math. Sci. Res. Inst. Publ., 60, Cambridge Univ. Press, Cambridge, 2013.
%

\bibitem{Hansen1} Hansen, P.~C., Discrete inverse problems: insight and algorithms, \emph{SIAM}, Philadelphia, 2010.

\bibitem{Hansen2} Hansen, P.~C., \emph{Analysis of discrete ill-posed problems by means of the L-curve}, SIAM review, 34(4), 561-580, 1992.

%\bibitem{H99} Helgason,~S. The Radon Transform, \emph{Birkh\"auser}, Basel, 1999.
%
%\bibitem{Helin}
%Helin,~T., Lassas,~M., Oksanen,~L., and Saksala,~T. \emph{Correlation based passive imaging with a white noise source.} Journal de Mathematiques Pures et Appliquees \textbf{116}(9), 132--160, 2018.
%
%
%
\bibitem{HUZ20} Hintz,~P., Uhlmann,~G. and Zhai,~J. \emph{An inverse boundary value problem for a semilinear wave equation on Lorentzian manifolds}, Int. Math. Res. Not., rnab088, 2021.

\bibitem{HUZ21} Hintz,~P., Uhlmann,~G. and Zhai,~J. \emph{The Dirichlet-to-Neumann map for a semilinear wave equation on Lorentzian manifolds}, arXiv preprint  	arXiv:2103.08110, 2021.


%\bibitem{IY01} Imanulov,~O. and Yamamoto,~M. \emph{Global uniqueness and stability in determining coefficients of wave equations}, Comm. PDE, \textbf{26}(7--8), 1409--1425, 2001.
%
%
%\bibitem{IsSun} Isakov,~V. and Sun,~Z. \emph{Stability estimates for hyperbolic inverse problems with local boundary data. } Inverse Problems \textbf{8}(2), 193--206, 1992.
%
%\bibitem{Isozaki} 
%Isozaki,~H., Kurylev,~L., and Lassas,~M. \emph{Conic singularities, generalized scattering matrix, and inverse scattering on asymptotically hyperbolic surfaces.} Journal fur die reine und angewandte Mathematik \textbf{724}, 53--103, 2017. 
%
%\bibitem{KKL01} Kachalov,~A.,  Kurylev,~Y., and Lassas,~M. Inverse boundary spectral problems. \emph{CRC Press}, 2001.
%
%\bibitem{KaNa02} Kang,~H. and Nakamura,~G. \emph{Identification of nonlinearity in a conductivity equation via the Dirichlet-to-Neumann map}, Inverse Problems, \textbf{18}(4), 1079, 2002. 
%
%
%\bibitem{KS14} Kenig,~C. and Salo,~M. \emph{Recent progress in the Cader\'on problem with partial data},  Contemp. Math., \textbf{615}, Amer. Math. Soc., Providence, RI, 2014.
%
\bibitem{KR2014} Knowles,~I. and Renka,~R.~J. \emph{Methods for numerical differentiation of noisy data}, Electron. J. Differ. Equ, 21, 235-246, 2014.

%\bibitem{KKLO}
%Kian,~Y., Kurylev,~Y., Lassas,~M., and Oksanen,~L. \emph{Unique recovery of lower order coefficients for hyperbolic equations from data on disjoint sets.} J. Differential Equations \textbf{267}(4), 2210--2238, 2019.
%
%
%\bibitem{KrKL}
%Krupchyk,~K., Kurylev,~Y., and Lassas,~M. {\em Inverse spectral problems on a closed manifold.} Journal de Mathematique Pures et Appliquees \textbf{90}, 42--59, 2008.
%
\bibitem{KrUh19} Krupchyk,~K. and Uhlmann,~G. \emph{Partial  data  inverse  problems  for  semilinear  elliptic  equations  with  gradient nonlinearities}, to appear in Mathematical Research Letters, arXiv:1909.08122, 2019
%
\bibitem{KrUh20} Krupchyk,~K. and Uhlmann,~G. \emph{A  remark  on  partial  data  inverse  problems  for  semilinear  elliptic  equations}, Proc. Amer. Math. Soc. \textbf{148}(2), 681--685, 2020.
%
%\bibitem{KLOU2014} Kurylev,~Y.,  Lassas,~M.,  Oksanen,~L.  and  Uhlmann,~G.  \emph{Inverse  problem  for  Einstein-scalar  field  equations}. arXiv:1406.4776, 2014
%
\bibitem{KLU18} Kurylev,~Y., Lassas,~M., and Uhlmann,~G.  \emph{Inverse problems for Lorentzian manifolds and non-linear hyperbolic equations.} Inventiones Mathematicae, \textbf{212}(3),781--857, 2018. \href{https://doi.org/10.1007/s00222-017-0780-y}{https://doi.org/10.1007/s00222-017-0780-y}
%
%
%\bibitem{KOP}
%Kurylev,~Y., Oksanen,~L., and Paternain,~G., \emph{Inverse problems for the connection Laplacian}, J. Differential Geom. \textbf{110}(3), 457--494, 2018. 
%
%\bibitem{LaUhYa20} Lai,~R-Y., Uhlmann,~G., and Yang,~Y. \emph{ Reconstruction of the collision kernel in the nonlinear boltzmann equation}, arXiv preprint arXiv:2003.09549, 2020.
%


\bibitem{LaUhYa20} Lai,~R-Y., Uhlmann,~G., and Yang,~Y. \emph{ Reconstruction of the collision kernel in the nonlinear boltzmann equation}, arXiv preprint arXiv:2003.09549, 2020.
%\bibitem{LLT86} Lasiecka,~I., Lions,~J-L., and Triggiani,~R. \emph{Non homogeneous boundary value problems for second order hyperbolic operators.} Journal de Mathématiques pures et Appliquées 65.2, 149--192, 1986.
%
%
%\bibitem{Lassas}
%Lassas,~M. \emph{Inverse problems for linear and non-linear hyperbolic equations.} Proc. Int. Congress of Math. ICM 2018, Rio de Janeiro, Brazil, Vol III, 3739--3760, 2018.
%
\bibitem{LLLS19a} Lassas,~M., Liimatainen,~T., Lin,~Y-H., and Salo,~M. \emph{Partial data inverse problems and simultaneous recovery of boundary and coefficients for semilinear elliptic equations}, arXiv preprint arXiv:1905.02764, 2019.
%
\bibitem{LLLS19b} Lassas,~M., Liimatainen,~T., Lin,~Y-H. and Salo,~M. \emph{Inverse problems for elliptic equations with power type nonlinearities}, arXiv preprint arXiv:1903.12562.

\bibitem{LLPT20} Lassas,~M, Liimatainen,~T., Potenciano-Machado,~L., and Tyni,~T. \emph{Uniqueness and stability of an inverse problem for a semi-linear wave equation}, arXiv preprint arXiv:2006.13193 (2020).

\bibitem{LLPT21} M. Lassas, T. Liimatainen, L. Potenciano-Machado and T. Tyni, \emph{Stability estimates for inverse problems for semi-linear wave equations on Lorentzian manifolds}, arXiv:2106.12257, (2021).

\bibitem{LSX21} L. Shuai, M. Salo and B. Xu, \textit{Increasing stability in the linearized inverse Schr\"{o} dinger potential problem with power type nonlinearities.}, arXiv preprint arXiv:2111.13446 (2021).

%
%\bibitem{LO}
%Lassas,~M. and Oksanen,~L., \emph{Inverse problem for the Riemannian wave equation with Dirichlet data and Neumann data on disjoint sets.} Duke Math. J. \textbf{163}, 1071--1103, 2014.
% 
\bibitem{LUW17}  Lassas,~M., Uhlmann,~G.  and Wang,~Y.  \emph{Determination of vacuum space-times from the Einstein-Maxwell equations}, arXiv preprint arXiv:1703.10704, 2017.
%
\bibitem{LUW18} Lassas,~M., Uhlmann,~G.  and Wang,~Y. \emph{Inverse Problems for Semilinear Wave Equations on Lorentzian Manifolds}, Comm. in Math. Phys., \textbf{360}, 555--609, 2018.
%
\bibitem{MG} Mitchell, A.R. and Griffiths~,~D.F., The Finite Difference Method in Partial Differential Equations, \emph{Wiley}, 1980.


%
\bibitem{MS2012}
Mueller J. and Siltanen S. {\em Linear and Nonlinear
  Inverse Problems with Practical Applications}, SIAM, Philadelphia, 2012.
%
%\bibitem{N96} Nachman,~A. \emph{Global uniqueness for a two-dimensional inverse boundary value problem}, Ann. of Math. \textbf{143}(2), 71--96, 1996.  
%
%
%\bibitem{NSU88} Nachman,~A., Sylvester,~J. and Uhlmann,~G. \emph{An n-dimensional Borg- Levinson theorem}, Comm. Math. Phys. \textbf{115}(4), 595--605, 1988.
%
%\bibitem{Na} Natterer,~F. The Mathematics of Computerized Tomography, \emph{Society for Industrial and Applied Mathematics (SIAM)}, Philadelphia, 2001.
%
% \bibitem{OSSU}
%Oksanen,~L., Salo,~M. Stefanov,~P., and Uhlmann,~G., \emph{Inverse problems for real principal type operators},arXiv preprint arXiv:2001.07599.
%
%\bibitem{Siltanen1}
%Siltanen,~S., Mueller,~J. and Isaacson,~D. \emph{An implementation of the
%  reconstruction algorithm of {A}. {N}achman for the 2-{D} inverse conductivity
%  problem}, Inverse Problems, \textbf{6}, 681--699, 2000.
%

\bibitem{Smith} Smith,~G.D. Numerical Solution of Partial Differential Equations: Finite Difference Methods, 3rd edition, \emph{Clarendon Press}, Oxford, 1985.

\bibitem{SSX21} L. Shuai, M. Salo and B. Xu, \textit{Increasing stability in the linearized inverse Schr\"{o} dinger potential problem with power type nonlinearities.}, arXiv preprint arXiv:2111.13446 (2021).


%\bibitem{Sogge} Sogge,~C.~D., Lectures on nonlinear wave equations, \emph{International Press Inc.}, 1995.
%
%\bibitem{St}
%Stefanov,~P. \emph{Inverse scattering problem for the wave equation with time-dependent potential}, J. Math. Anal. Appl., \textbf{140}, 351--362, 1989.
%
%\bibitem{StUh} Stefanov,~P. and Uhlmann,~G. \emph{Stable Determination of Generic Simple Metrics from the Hyperbolic Dirichlet-to-Neumann Map}, Int. Math. Res. Not., \textbf{17}, 1047--1061, 2005.
%
%\bibitem{SY}
%Stefanov,~P. and Yang,~Y. \emph{The inverse problem for the Dirichlet-to-Neumann map on Lorentzian manifolds}. Anal. PDE., \textbf{11}(6), 1381--1414, 2018.
%
%\bibitem{SunUh97} Sun,~Z. and Uhlmann,~G. \emph{Inverse problems in quasilinear anisotropic media}, American J. Math. \textbf{119}(4), 771--797, 1997.
%
%\bibitem{SyUh87} Sylvester,~J. and Uhlmann,~G. \emph{A Global Uniqueness Theorem for an Inverse Boundary Value Problem}, Annals of Math. \textbf{125}(1), 153--169, 1987. 
%
%\bibitem{Uh13} Uhlmann,~G. \emph{30 Years of Calder\'on's Problem},  S\'emin. \'Equ. D\'eriv. Partielles, \'Ecole Polytech., Palaiseau, 2014. 
%
%\bibitem{UhWa18} Uhlmann,~G. and Wang,~Y. \emph{Determination of space-time structures from gravitational perturbations}, To appear in Comm. Pure App. Math. (2018).
%

\bibitem{UZ21} Uhlmann, G. and Zhai, J., \textit{On an inverse boundary value problem for a nonlinear elastic wave equation}, J. Math. Pures Appl., \textbf{153}, (2021), 114--136, \url{https://doi.org/10.1016/j.matpur.2021.07.005}.

\bibitem{UZ21b} Uhlmann, G. and Zhai, J., \textit{Inverse problems for nonlinear hyperbolic equations}, Discrete Contin. Dyn. Syst., \textbf{41}, (2021), No. 1, 455--469, \url{https://doi.org/10.3934/dcds.2020380}.

%\bibitem{Ze95} Zeidler,~E. Applied functional analysis, \emph{Springer-Verlag}, New York, 1995.

\bibitem{WZ2019} Wang,~Y. and  Zhou,~T. \emph{Inverse problems for quadratic derivative nonlinear wave equations}. Comm. PDE, \textbf{44}(11), 1140--1158, 2019.

\end{thebibliography}
\end{document}